\documentclass{article}
\usepackage{amsmath,amssymb,enumerate,bbm}

\catcode`\à=\active \defà{\`a} \catcode`\À=\active \defÀ{\`A}
\catcode`\á=\active \defá{\'a} \catcode`\Á=\active \defÁ{\'A}
\catcode`\ä=\active \defä{\"a} \catcode`\Ä=\active \defÄ{\"A}
\catcode`\â=\active \defâ{\^a} \catcode`\Â=\active \defÂ{\^A}
\catcode`\å=\active \defå{{\aa}} \catcode`\Å=\active \defÅ{{\AA}}
\catcode`\ç=\active \defç{\c{c}} \catcode`\Ç=\active \defÇ{\c{C}}
\catcode`\è=\active \defè{\`e} \catcode`\È=\active \defÈ{\`E}
\catcode`\é=\active \defé{\'e} \catcode`\É=\active \defÉ{\'E}
\catcode`\ë=\active \defë{\"e} \catcode`\Ë=\active \defË{\"E}
\catcode`\ê=\active \defê{\^e} \catcode`\Ê=\active \defÊ{\^E}
\catcode`\ì=\active \defì{\`{\i}} \catcode`\Ì=\active \defÌ{\`{\I}}
\catcode`\í=\active \defí{\'{\i}} \catcode`\Í=\active \defÍ{\'{\I}}
\catcode`\ï=\active \defï{\"{\i}} \catcode`\Ï=\active \defÏ{\"{\I}}
\catcode`\î=\active \defî{\^{\i}} \catcode`\Î=\active \defÎ{\^{\I}}
\catcode`\ò=\active \defò{\`o} \catcode`\Ò=\active \defÒ{\`O}
\catcode`\ó=\active \defó{\'o} \catcode`\Ó=\active \defÓ{\'O}
\catcode`\ö=\active \defö{\"o} \catcode`\Ö=\active \defÖ{\"O}
\catcode`\ô=\active \defô{\^o} \catcode`\Ô=\active \defÔ{\^O}
\catcode`\ù=\active \defù{\`u} \catcode`\Ù=\active \defÙ{\`U}
\catcode`\ú=\active \defú{\'u} \catcode`\Ú=\active \defÚ{\'U}
\catcode`\ü=\active \defü{\"u} \catcode`\Ü=\active \defÜ{\"U}
\catcode`\û=\active \defû{\^u} \catcode`\Û=\active \defÛ{\^U}
\catcode`\ý=\active \defý{\'y} \catcode`\Ý=\active \defÝ{\'Y}
\catcode`\ÿ=\active \defÿ{\"y} \catcode`\˜=\active \def˜{\"Y}
\catcode`\½=\active \def½{!`}
\catcode`\¾=\active \def¾{?`}
\catcode`\ß=\active \defß{{\ss}}
\newcommand{\tmop}[1]{\ensuremath{\operatorname{#1}}}
\newtheorem{theorem}{Theorem}
\newtheorem{definition}{Definition}
\newtheorem{proposition}{Proposition}
\newenvironment{proof}{\noindent\textbf{Proof\ }}{\hspace*{\fill}$\Box$\medskip}
\newcommand{\tmmathbf}[1]{\ensuremath{\boldsymbol{#1}}}
\newcommand{\tmscript}[1]{\text{\scriptsize $#1$}}
\newcommand{\di}{\displaystyle}

\begin{document}

\title{Formulas for the Connes-Moscovici Hopf algebra} 
\author{Frédéric Menous} 
\date{}
\maketitle

\begin{abstract}
  We give explicit formulas for the coproduct and the antipode in the
  Connes-Moscovici Hopf algebra $\mathcal{H}_{\tmop{CM}}$. To do so, we first
  restrict ourselves to a sub-Hopf algebra $\mathcal{H}^1_{\tmop{CM}}$
  containing the nontrivial elements, namely those for which the coproduct and
  the antipode are nontrivial. There are two ways to obtain explicit
  formulas. On 
  one hand, the algebra $\mathcal{H}^1_{\tmop{CM}}$ is isomorphic to the Faà
  di Bruno Hopf algebra of coordinates on the group of identity-tangent
  diffeomorphism and computations become easy using substitution automorphisms
  rather than diffeomorphisms. On the other hand, the algebra
  $\mathcal{H}^1_{\tmop{CM}}$ is isomorphic to a sub-Hopf algebra of the
  classical shuffle Hopf algebra which appears naturally in resummation
  theory, in the framework of formal and analytic  conjugacy of vector fields.
  Using the very simple structure of the shuffle Hopf algebra, we derive once
  again explicit formulas for the coproduct and the antipode in 
  $\mathcal{H}^1_{\tmop{CM}}$.
\end{abstract}

\section{Introduction.}

The Connes-Moscovici Hopf algebra $\mathcal{H}_{\tmop{CM}}$ was introduced in
{\cite{cm}} in the context of noncommutative geometry. Because of its relation
with the Lie algebra of formal vector fields, it was also proved in
{\cite{cm}} that its subalgebra $\mathcal{H}^1_{\tmop{CM}}$ is isomorphic to
the Faà di Bruno Hopf algebra of coordinates of identity-tangent
diffeomorphisms (see {\cite{cm}},{\cite{fig}}). In the past years, it appeared
that this Hopf algebra was strongly related to the Hopf algebras of trees (see
{\cite{ck0}}) or graphs (see {\cite{ck1}},{\cite{ck2}}) underlying
perturbative renormalization in quantum field theory.

Our aim is to give explicit formulas for the coproduct and the antipode in
$\mathcal{H}^1_{\tmop{CM}}$, since only recursive formulas seem to be known.

We remind in section \ref{cm} the definition of the Connes-Moscovici Hopf
algebra, as well as its properties and links with the Faà di Bruno Hopf
algebra and identity-tangent diffeomorphisms (for details, see
{\cite{cm}},{\cite{fig}}). The formulas are given in section \ref{for}. We
present a proof based on the isomorphism between identity-tangent
diffeomorphisms and substitution automorphisms which are easier to handle in
the computations. These manipulations on substitution automorphisms are very
common in J. Ecalle's work on the formal classification of differential
equations, vector fields, diffeomorphism... (see
{\cite{et1}},{\cite{et2}},{\cite{et3}},{\cite{esn}}). In fact, the first proof
for these formulas was based on mould calculus and shuffle Hopf algebras,
which we shortly describe in section \ref{moulds}. Sections \ref{mor} and
\ref{ini} give the outlines of the initial proof based on a Hopf morphism from
$\mathcal{H}^1 \subset \mathcal{H}_{\tmop{CM}}$ in a shuffle Hopf algebra. 

{\tableofcontents}

\section{Connes-Moscovici and Faà di Bruno Hopf algebras.\label{cm}}

\subsection{The Connes-Moscovici Hopf algebra}

The Connes-Moscovici Hopf algebra $\mathcal{H}_{\tmop{CM}}$ defined in
{\cite{cm}} is the enveloping algebra of the Lie algebra which is the linear
span of $Y$, $X$, $\delta_n$, $n \geq 1$ with the relations,
\begin{equation}
  [ X, Y ] = X, [ Y, \delta_n ] = n \delta_n, [ \delta_n, \delta_m ] = 0, [ X,
  \delta_n ] = \delta_{n + 1} \label{crocm}
\end{equation}
for all $m, n \geq 1$. The coproduct $\Delta$ in $\mathcal{H}_{\tmop{CM}}$ is
defined by
\begin{equation}
  \Delta ( Y ) = Y \otimes 1 + 1 \otimes Y, \Delta ( X ) = X \otimes 1 + 1
  \otimes X + \delta_1 \otimes Y, \Delta ( \delta_1 ) = \delta_1 \otimes 1 + 1
  \otimes \delta_1
\end{equation}
where $\Delta ( \delta_n )$ is defined recursively, using equation \ref{crocm}
and the identity
\begin{equation}
  \forall h_1, h_2 \in \mathcal{H}_{\tmop{CM}}, \quad \Delta ( h_1 h_2 ) =
  \Delta ( h_1 ) \Delta ( h_2 )
\end{equation}
The coproduct of $X$ and $Y$ is given, whereas the coproduct of $\delta_n$ is
nontrivial. Nonetheless, the algebra generated by $\{ \delta_n, \hspace{1em} n
\geq 1 \}$ is a graded sub-Hopf algebra $\mathcal{H}^1_{\tmop{CM}} \subset
\mathcal{H}_{\tmop{CM}}$ where the graduation is defined by
\begin{equation}
  \tmop{gr} ( \delta_{n_1} \ldots \delta_{n_s} ) = n_1 + \ldots + n_s
\end{equation}
As mentioned in {\cite{cm}}, the Hopf algebra $\mathcal{H}^1_{\tmop{CM}}$ is
strongly linked to Faà di Bruno Hopf algebra.

\subsection{The Faà di Bruno Hopf algebra}

Let us consider the group of formal identity tangent diffeomorphisms :
\[ G_2 = \{ f ( x ) = x + \sum_{n \geq 1} f_n x^{n + 1} \in \mathbbm{R}[ [ x ]
   ] \} \]
with, by convention, the product $\mu : G_2 \times G_2 \rightarrow G_2$ :
\[ \mu ( f, g ) = g \circ f \]
For $n \geq 0$, the functionals on $G_2$ defined by
\[ a_n ( f ) = \frac{1}{( n + 1 ) !_{}} ( \partial_x^{n + 1} f ) ( 0 ) = f_n
   \hspace{1em} a_n : G_2 \rightarrow \mathbbm{R} \]
are called de Faà di Bruno coordinates on the group $G_2$ and $a_0 = 1$ being
the unit, they generates a graded unital commutative algebra
\[ \mathcal{H}_{\tmop{FdB}} =\mathbbm{R}[ a_1, \ldots, a_n, \ldots ]
   \hspace{1em} ( \tmop{gr} ( a_n ) = n ) \]
Moreover, the action of these functionals on a product in $G_2$ defines a
coproduct on $\mathcal{H}_{\tmop{FdB}}$ that turns to be a graded connected
Hopf algebra (see {\cite{fig}} for details). For $n \geq 0$, the coproduct is
defined by
\begin{equation}
  a_n \circ \mu = m \circ \Delta ( a_n ) \label{copfdb}
\end{equation}
where $m$ is the usual multiplication in $\mathbbm{R}$, and the antipode reads
\[ S \circ a_n = a_n \circ \tmop{rec} \]
where $\tmop{rec} ( \varphi ) = \varphi^{- 1}$ is the composition inverse of
$\varphi$.

For example if $f ( x ) = x + \sum_{n \geq 1} f_n x^{n + 1}$ and $g ( x ) = x
+ \sum_{n \geq 1} g_n x^{n + 1}$ then if $h = \mu ( f, g ) = g \circ f$ and $h
( x ) = x + \sum_{n \geq 1} h_n x^{n + 1}$,
\[ \begin{array}{ccccccc}
     a_0 ( h ) & = & 1 = a_0 ( f ) a_0 ( g ) & \rightarrow & \Delta a_0 & = &
     a_0 \otimes a_0\\
     a_1 ( h ) & = & f_1 + h_1 & \rightarrow & \Delta a_1 & = & a_1 \otimes
     a_0 + a_0 \otimes a_1\\
     a_2 ( h ) & = & f_2 + f_1 g_1 + g_2 & \rightarrow & \Delta a_2 & = & a_2
     \otimes a_0 + a_1 \otimes a_1 + a_0 \otimes a_2
   \end{array} \]
As proved in {\cite{cm}} and {\cite{fig}}, there exists a Hopf isomorphism
between $\mathcal{H}_{\tmop{FdB}}$ and $\mathcal{H}^1_{\tmop{CM}}$.

\subsection{Connes-Moscovici coordinates}

Following {\cite{cm}}, one can define new functionals on $G_2 $ by $\gamma_0 =
a_0$=1 (unit) and for $n \geq 1$,
\[ \gamma_n ( f ) = ( \partial^n_x \log ( f' ) ) ( 0 ) \]
These functionals, which may be called the Connes-Moscovici coordinates on
$G_2$, freely generates the Faà di Bruno Hopf algebra :
\[ \mathcal{H}_{\tmop{FdB}} =\mathbbm{R}[ a_1, \ldots, a_n, \ldots ]
   =\mathbbm{R}[ \gamma_1, \ldots, \gamma_n, \ldots ] \hspace{1em} \tmop{gr} (
   a_n ) = \tmop{gr} ( \gamma_n ) = n \]
and their coproduct is given by the formula \ref{copfdb}. Now, see
{\cite{cm}}, {\cite{ck0}} :

\begin{theorem}
  The map $\Theta$ defined by $\Theta ( \delta_n ) = \gamma_n$ is a graded
  Hopf isomorphism between $\mathcal{H}_{\tmop{FdB}}$ and
  $\mathcal{H}^1_{\tmop{CM}}$
\end{theorem}

This means that the coproduct and the antipode in $\mathcal{H}^1_{\tmop{CM}}$
can be rather computed in $\mathcal{H}_{\tmop{FdB}}$. Unfortunately, if the
coproduct and the antipode is well-known for the functionals $a_n$, using the
Faà di Bruno formulas for the composition and the inverse of diffeomorphisms
in $G_2$, it seems that formulas for the $\gamma_n$ cannot be easily derived.
In order to do so, we will either work with substitution automorphism which
are easier to handle than diffeomorphisms (see section \ref{for}, or identify
$\mathcal{H}_{\tmop{FdB}}$ as a sub-Hopf algebra of a shuffle Hopf algebra and
use mould calculus (see sections \ref{moulds}, \ref{mor}, \ref{ini}).

\section{Formulas in $\mathcal{H}^1_{\tmop{CM}}$.\label{for}}

\subsection{Notations}

In the sequel we note
\[ \mathcal{N}= \{ \tmmathbf{n}= ( n_1, \ldots, n_s ) \in (\mathbbm{N}^{\ast}
   )^s, \hspace{1em} s \geq 1 \} \]
For $\tmmathbf{n}= ( n_1, \ldots, n_s ) \in \mathcal{N}$,
\[ \| \tmmathbf{n} \| = n_1 + \ldots + n_s, \hspace{1em} l ( \tmmathbf{n} ) =
   s \]
and if $n \geq 1$,
\[ \mathcal{N}_n = \{ \tmmathbf{n} \in \mathcal{N} \hspace{1em} ; \hspace{1em}
   \| \tmmathbf{n} \| = n \} \]
For a tuple $\tmmathbf{n}= ( n_1, \ldots, n_s ) \in \mathcal{N}$, we note
$\tmmathbf{n}! = n_1 ! \ldots n_s !$. More over, $\tmop{Split} (\tmmathbf{n})$
is the subset of $\bigcup_{t \geq 1} \mathcal{N}^t$ such that $(\tmmathbf{n}^1,
\ldots,\tmmathbf{n}^t ) \in \tmop{Split} (\tmmathbf{n})$ if and only if the
concatenation of $(\tmmathbf{n}^1, \ldots,\tmmathbf{n}^t )$ is equal to
$\tmmathbf{n}$ :
\begin{equation}
  \tmop{Split} (\tmmathbf{n}) = \{ (\tmmathbf{n}^1, \ldots,\tmmathbf{n}^t )
  \in \mathcal{N}^t, \hspace{1em} \tmmathbf{n}^1 \ldots \tmmathbf{n}^t
  =\tmmathbf{n} \}
\end{equation}

In summation formulas, we will use the fact that
\begin{equation}
  \bigcup_{\tmmathbf{n} \in \mathcal{N}_n} \tmop{Split} (\tmmathbf{n}) =
  \bigcup_{\tmmathbf{n}= ( n_1, \ldots, n_s ) \in \mathcal{N}_n}
  \mathcal{N}_{n_1} \times \ldots \times \mathcal{N}_{n_s}
\end{equation}
so that if $f$ is a function on $\mathcal{N}$ and $g$ is a function on
$\bigcup_{t \geq 1} \mathcal{N}^t$, for $n \geq 1$,
\begin{equation}
  \begin{array}{l}
   \displaystyle \sum_{\tmmathbf{n}= ( n_1, \ldots, n_s ) \in \mathcal{N}_n}
    \sum_{\tmscript{\begin{array}{c}
      \tmmathbf{m}^1 \in \mathcal{N}_{n_1}\\
      { \vdots}\\
     \tmmathbf{m}^s \in \mathcal{N}_{n_s}
    \end{array}}} f (\tmmathbf{n}) g ( \tmmathbf{m}^1, \ldots, \tmmathbf{m}^s
    ) =\\
    \hspace{10em} \displaystyle \sum_{\tmmathbf{n} \in \mathcal{N}_n}
    \sum_{\tmmathbf{m}^1 
    \ldots \tmmathbf{m}^s =\tmmathbf{n}} f ( \| \tmmathbf{m}^1 \|, \ldots, \|
    \tmmathbf{m}^s \| ) g ( \tmmathbf{m}^1, \ldots, \tmmathbf{m}^s )
  \end{array}
\end{equation}
where $\displaystyle \sum_{\tmmathbf{m}^1 \ldots \tmmathbf{m}^s
  =\tmmathbf{n}}$is the sum 
over $\tmop{Split} (\tmmathbf{n})$.

Finally, for $(\tmmathbf{n}^1, \ldots,\tmmathbf{n}^t ) \in \mathcal{N}^t$ ($t
\geq 1$),
\begin{equation}
  A (\tmmathbf{n}^1, \ldots,\tmmathbf{n}^t ) = \frac{1}{l (\tmmathbf{n}^1 ) !
  \ldots l (\tmmathbf{n}^t ) !} \prod_{i = 1}^t \frac{1}{\| \tmmathbf{n}^i \|
  + 1}
\end{equation}
and, for $k \geq 1$,
\begin{equation}
  B_k ( \tmmathbf{n}^1, \ldots,^{} \tmmathbf{n}^t ) = C^{l ( \tmmathbf{n}^t
  )}_k \prod_{i = 1}^{t - 1} C_{\| \tmmathbf{n}^{i + 1} \| + \ldots + \|
  \tmmathbf{n}^t \| + k}^{l ( \tmmathbf{n}^i )}
\end{equation}

\subsection{Main formulas}

We will now prove the following formulas :

\begin{theorem}
  \label{th2}For $n \geq 1$,
  \begin{equation}
\begin{array}{rcl} 
   \Delta ( \delta_n )& =& \delta_n \otimes 1 + 1 \otimes \delta_n \\
 & &\displaystyle+    \sum_{\tmscript{\begin{array}{c}
      ( n_1, \ldots, n_{s + 1} ) \in \mathcal{N}_n\\
      s \geq 1
    \end{array}}} \frac{n!}{n_1 ! \ldots n_{s + 1} !} \alpha^{n_1, \ldots,
    n_s}_{n_{s + 1}}  \delta_{n_1} \ldots \delta_{n_s}  \otimes \delta_{n_{s
    + 1}}\end{array}
  \end{equation}
  and, for $\tmmathbf{n}= ( n_1, \ldots, n_s ) \in \mathcal{N}$ ($l
  (\tmmathbf{n}) = s$) and $m \geq 1$,
  \begin{equation}
    \alpha^{\tmmathbf{n}}_m = \sum_{t = 1}^{l (\tmmathbf{n})} C^{t_{}}_m
    \sum_{\tmscript{\begin{array}{c}
      \tmmathbf{n}^1 \ldots \tmmathbf{n}^t =\tmmathbf{n}
    \end{array}}} A (\tmmathbf{n}^1, \ldots,\tmmathbf{n}^t )
  \end{equation}
  where, for $\tmmathbf{n}= ( n_1, \ldots, n_s ) \in \mathcal{N}$, $l
  (\tmmathbf{n}) = 1$, $\| \tmmathbf{n} \| = n_1 + \ldots + n_s$ and with the
  convention $C^t_m = \displaystyle \frac{m!}{t! ( m - t ) !} = 0$ if $t > m$.
\end{theorem}

For the antipode $S$ :

\begin{theorem}
  \label{th3}For $n \geq 1$,
  \begin{equation}
    S ( \delta_n ) = \sum_{\tmscript{\begin{array}{c}
      \tmmathbf{n}= ( n_1, \ldots, n_s ) \in \mathcal{N}\\
      n_1 + \ldots + n_s = n
    \end{array}}} \frac{n!}{n_1 ! \ldots n_s !} \beta^{n_1, \ldots, n_s}
    \delta_{n_1} \ldots \delta_{n_s}
  \end{equation}
  with $\beta^{n_1} = - 1$ and, if $\tmmathbf{n}= ( n_1, \ldots, n_{s + 1} )
  \in \mathcal{N}$ ($s \geq 1$),
  \begin{equation}
    \beta^{n_1, \ldots, n_s, n_{s + 1}} = \sum_{t = 1}^s
    \sum_{\tmscript{\begin{array}{c}
      \tmmathbf{n}^1 \ldots \tmmathbf{n}^t =\tmmathbf{n}
    \end{array}}} U^{\| \tmmathbf{n}^1 \|, \ldots, \| \tmmathbf{n}^t \|}_{n_{s
    + 1}} A (\tmmathbf{n}^1, \ldots,\tmmathbf{n}^t )
  \end{equation}
  where, if $\tmmathbf{m}= ( m_1, \ldots, m_t ) \in \mathcal{N}/ \{ \emptyset
  \}$ and $k \geq 1$,
  \begin{equation}
    U^{\tmmathbf{m}}_k = \sum_{i = 1}^{l (\tmmathbf{m})} ( - 1 )^{i - 1}
    \sum_{\tmscript{\begin{array}{c}
      \tmmathbf{m}^1 \ldots \tmmathbf{m}^i =\tmmathbf{m}
    \end{array}}} B_k ( \tmmathbf{m}^1, \ldots,^{} \tmmathbf{m}^i )
  \end{equation}
\end{theorem}

We will now give the more recent proof of this formulas. These formulas were
first conjectured and then proved using a Hopf morphism between
$\mathcal{H}^1_{\tmop{CM}}$ and a shuffle Hopf algebra noted $\tmop{sh}
(\mathbbm{N}^{\ast} )$. We will come back later on this morphism and the
afferent proofs. Let us first look at the correspondence between FdB
coordinates and the CM coordinates on $G_2$.

\subsection{Coordinates on $G_2$}

Let $\varphi ( x ) = x +\displaystyle \sum_{n \geq 1} \varphi_n x^{n +
  1}$. We have for $n
\geq 1$ :
\begin{equation}
  a_n ( \varphi ) = \varphi_n, \hspace{1em} \gamma_n ( \varphi ) = (
  \partial^n_x \log ( \varphi' ) ) ( 0 ) = f_n
\end{equation}
If $f ( x ) = \displaystyle \sum_{n \geq 1} \frac{f_n}{n!} x^n$, then
\begin{equation}
  f ( x ) = \log ( \varphi' ( x ) ) \hspace{1em} \varphi ( x ) = \int_0^x e^{f
  ( t )} d t \label{corr}
\end{equation}
For any sequence $( u_n )_{n \geq 1}$, we note
\begin{equation}
  \forall \tmmathbf{n}= ( n_1, \ldots, n_s ) \in \mathcal{N}, \hspace{1em}
  u_{\tmmathbf{n}} = u_{n_1} \ldots u_{n_s}
\end{equation}
Using equation \ref{corr}, we get easily that
\begin{equation}
  \begin{array}{ccc}
    f ( x ) & = &\displaystyle  \sum_{\tmmathbf{n}= ( n_1, \ldots, n_s ) \in \mathcal{N}}
    \frac{( - 1 )^{l ( \tmmathbf{n} )}}{l ( \tmmathbf{n} )} ( n_1 + 1 ) \ldots
    ( n_s + 1 ) \varphi_{\tmmathbf{n}} x^{\| \tmmathbf{n} \|}\\
    \varphi ( x ) & = & \displaystyle x + \sum_{\tmmathbf{n}= ( n_1, \ldots, n_s ) \in
    \mathcal{N}} \frac{1}{l ( \tmmathbf{n} ) ! \tmmathbf{n} !} 
    \frac{f_{\tmmathbf{n}}}{\| \tmmathbf{n} \| + 1} x^{\| \tmmathbf{n} \| + 1}
  \end{array} \label{e19}
\end{equation}
and these formulas establish the correspondence between FdB and CM coordinates
on $G_2$. In order to prove theorems \ref{th2} and \ref{th3}, we need to
understand how these coordinates read on $\varphi^{- 1}$ and $\mu ( \varphi,
\psi ) = \psi \circ \varphi$ ($\varphi, \psi \in G_2$). To do so, we will
rather work with substitution automorphisms than with diffeomorphism.

\subsection{Taylor expansions and substitution automorphisms\label{aut}}

\begin{definition}
  Let $\tilde{G}_2$ be the set of linear maps from $\mathbbm{R}[ [ x ] ]$ to
  $\mathbbm{R}[ [ x ] ]$ such that
  \begin{itemize}
    \item[1.] For $F \in \tilde{G}_2$, the image $F ( x )$ by $F$ of the series
    $x$ is in $G_2$.
    
    \item[2.] For any two series $A$ and $B$ in $\mathbbm{R}[ [ x ] ]$, we have
    \begin{equation}
      F ( A.B ) = F ( A ) .F ( B ) \label{prod}
    \end{equation}
  \end{itemize}
\end{definition}

The elements of $\tilde{G}_2$ are called substitution automorphisms and

\begin{theorem}
  $\tilde{G}_2$ is a group for the composition and the map :
  \[ \begin{array}{ccccc}
       \tau & : & \text{$\tilde{G}_2$} & \rightarrow & G_2\\
       &  & F & \mapsto & \varphi ( x ) = F ( x )
     \end{array} \]
  defines an isomorphism between the groups $\tilde{G}_2$ and $G_2$. Moreover,
  for $A \in \mathbbm{R}[ [ x ] ]$,
  \begin{equation}
    F ( A ) = A \circ \tau ( F )
  \end{equation}
\end{theorem}

\begin{proof}
  If $F \in \widetilde{G_{}}_2$, then, thanks to equation \ref{prod}, for $k
  \geq 0$,
  \begin{equation}
    \text{$F ( x^k ) = \left( F ( x ) \right)^k = ( \tau ( F ) ( x ) )^k = (
    \varphi ( x ) )^k$}
  \end{equation}
  thus, for $A ( x ) = \sum_{k \geq 0} A_k x^k \in \mathbbm{R} [ [ x ] ]$,
  \begin{equation}
    \begin{array}{ccc}
      F ( A ) ( x ) & = & F \left( \sum_{k \geq 0} A_k x^k \right)\\
      & = & \sum_{k \geq 0} A_k F ( x^k )\\
      & = & \sum_{k \geq 0} A_k ( \varphi ( x ) )^k\\
      & = & A \circ ( \tau ( F ) ) ( x )
    \end{array}
  \end{equation}
  This proves that $\tau$ is injective and for any $\varphi \in G_2$ the map
  \[ \begin{array}{ccccc}
       F & : & \mathbbm{R}[ [ x ] ] & \rightarrow & \mathbbm{R}[ [ x ] ]\\
       &  & A & \mapsto & A \circ \varphi
     \end{array} \]
  is a substitution automorphism of $\tilde{G}_2$ such that $\tau ( F ) =
  \varphi$. The map $\tau$ is a bijection. Now, for $F$ and $G$ in
  $\tilde{G}_2$,
  \begin{equation}
    \tau ( F \circ G ) ( x ) = F ( G ( x ) ) = \tau ( G ) \circ \tau ( F ) ( x
    ) = \mu ( \tau ( F ), \tau ( G ) ) ( x )
  \end{equation}
  and if $H = \tau^{- 1} ( ( \tau ( F ) )^{- 1} )$ then $F \circ H = H \circ F
  = \tmop{Id}$. This ends the proof.
  
  $^{}$
\end{proof}

Using Taylor expansion, we also get formulas for $\tau^{- 1} ( \varphi )$,
$\varphi \in G_2$,

\begin{proposition}
  Let $\varphi ( x ) = \di x + \sum_{n \geq 1} \varphi_n x^{n + 1} \in G_2$ and $F
  = \tau^{- 1} ( \varphi )$, then
  \begin{equation}
    F = \tmop{Id} + \sum_{\tmmathbf{n}= ( n_1, \ldots, n_s ) \in \mathcal{N}}
    \frac{1}{l ( \tmmathbf{n} ) !} \varphi_{\tmmathbf{n}} x^{\| \tmmathbf{n}
    \| + l ( \tmmathbf{n} )} \partial_x^{l ( \tmmathbf{n} )}
  \end{equation}
\end{proposition}

This also means that $F$ can be decomposed in homogeneous components :
\begin{equation}
  F = \tmop{Id} + \sum_{n \geq 1} F_n \hspace{1em}, \hspace{1em} F_n =
  \sum_{\tmmathbf{n}= ( n_1, \ldots, n_s ) \in \mathcal{N}_n} \frac{1}{l (
  \tmmathbf{n} ) !} \varphi_{\tmmathbf{n}} x^{\| \tmmathbf{n} \| + l (
  \tmmathbf{n} )} \partial_x^{l ( \tmmathbf{n} )}
\end{equation}
such that
\begin{equation}
  \forall n \geq 1, \hspace{1em} \forall k \geq 1, \hspace{1em} \exists c \in
  \mathbbm{R}, \hspace{1em} F_n ( x^k ) = c x^{n + k}
\end{equation}
\begin{proof}
  If $\varphi ( x ) = x + \di \sum_{n \geq 1} \varphi_n x^{n + 1} = x +
  \bar{\varphi} ( x ) \in G_2$, then, if $F = \tau^{- 1} ( \varphi )$, then
  for $A \in \mathbbm{R} [ [ x ] ]$,
  \[ \begin{array}{ccc}
       F ( A ) ( x ) & = & A ( x + \bar{\varphi} ( x ) )\\
       & = & \di A ( x ) + \sum_{s \geq 1} \frac{\left( \bar{\varphi} ( x ) )^s
       \right.}{s!} A^{( s )} ( x )\\
       & = & \di A ( x ) + \sum_{s \geq 1} \sum_{n_1 \geq 1, \ldots, n_s \geq 1}
       \frac{1}{s!} \varphi_{n_1} \ldots \varphi_{n_s} x^{n_1 + \ldots + n_s +
       s} A^{( s )} ( x )\\
       & = & \di \left( \tmop{Id} + \sum_{\tmmathbf{n}= ( n_1, \ldots, n_s ) \in
       \mathcal{N}} \frac{1}{l ( \tmmathbf{n} ) !} \varphi_{\tmmathbf{n}}
       x^{\| \tmmathbf{n} \| + l ( \tmmathbf{n} )} \partial_x^{l (
       \tmmathbf{n} )} \right) \left( A ( x ) \right)
     \end{array} \]
  
\end{proof}

The automorphism $F$ can be seen as a differential operator acting on
$\mathbbm{R}[ [ x ] ]$ and from now on we note multiplicatively the action of
such operators :
\begin{equation}
  F. \varphi = F ( \varphi )
\end{equation}

As this will be of some use later, let us give the following formula : If
$\tmmathbf{n}= ( n_1, \ldots, n_s ) \in \mathcal{N}_{}$ and $k \geq 1$,
\begin{equation}
  \begin{array}{ccc}
    F_{\tmmathbf{n}} .x^k & = & F_{n_1} \ldots F_{n_s} .x^k\\
    & = & \di \sum_{\tmscript{\begin{array}{c}
      \tmmathbf{m}^i \in \mathcal{N}_{n_i}\\
      1\leq i\leq s\\
      \end{array}}} \left( \frac{\varphi_{\tmmathbf{m}^1}x^{n_1 + l (
        \tmmathbf{m}^1 )}}{ l ( \tmmathbf{m}^1
      ) !}  \partial_x^{l (
    \tmmathbf{m}^1 )} \right) \ldots \left(
  \frac{\varphi_{\tmmathbf{m}^s}x^{n_s + l (
      \tmmathbf{m}^s )} }{l
    ( \tmmathbf{m}^s ) !} 
    \partial_x^{l ( \tmmathbf{m}^s )} \right) .x^k\\
    & = & \di \sum_{\tmscript{\begin{array}{c}
      \tmmathbf{m}^i \in \mathcal{N}_{n_i}\\
      1\leq i \leq s\\
      \end{array}}} B_k ( \tmmathbf{m}^1, \ldots,^{} \tmmathbf{m}^s )
    \varphi_{\tmmathbf{m}^1} \ldots \varphi_{\tmmathbf{m}^s} x^{\|
    \tmmathbf{n} \| + k}
  \end{array}
\end{equation}
where
\[ B_k ( \tmmathbf{m}^1, \ldots,^{} \tmmathbf{m}^s ) = C^{l ( \tmmathbf{m}^s
   )}_k \prod_{i = 1}^{s - 1} C_{\| \tmmathbf{m}^{i + 1} \| + \ldots + \|
   \tmmathbf{m}^s \| + k}^{l ( \tmmathbf{m}^i )} \]
With these results one can already derive formulas for the $\tmop{FdB}$
coordinates on $G_2$.

\subsection{Formulas in $\mathcal{H}_{\tmop{FdB}}$}

We recover the usual formulas :

\begin{proposition}
  \label{prop7}We have for $n \geq 1$,
  \begin{equation}
    \Delta ( a_n ) = a_n \otimes 1 + 1 \otimes a_n + \sum_{k = 1}^{n - 1}
    \sum_{\tmmathbf{n}= ( n_1, \ldots, n_s ) \in \mathcal{N}_k} C^{l (
    \tmmathbf{n} )}_{n - k + 1} a_{\tmmathbf{n}} \otimes a_{n - k}
  \end{equation}
  and
  \begin{equation}
    S ( a_n ) = \sum_{\tmmathbf{n} \in \mathcal{N}_n} \left(
    \sum_{\tmmathbf{m}^1 \ldots \tmmathbf{m}^s = \tmmathbf{n}} ( - 1 )^s B_1 (
    \tmmathbf{m}^1, \ldots,^{} \tmmathbf{m}^s ) \right) a_{\tmmathbf{n}}
  \end{equation}
  
\end{proposition}

\begin{proof}
  Let $\varphi ( x ) = x + \di \sum_{n \geq 1} \varphi_n x^{n + 1}$ and $\psi ( x
  ) = x + \di \sum_{n \geq 1} \psi_n x^{n + 1}$ two elements of $G_2$ and $\eta =
  \mu ( \varphi, \psi ) = \psi \circ \varphi$ with
  \begin{equation}
    \eta ( x ) = x + \sum_{n \geq 1} \eta_n x^{n + 1}
  \end{equation}
  If $F$, $G$ and $H$ are the substitution automorphisms corresponding to
  $\varphi$, $\psi$ and $\eta$, then $H = F \circ G$ :
  \begin{equation}
    \begin{array}{cccc}
      H & = & \di \tmop{Id} + \sum_{n \geq 1} H_n & \\
      & = & \di \left( \tmop{Id} + \sum_{n \geq 1} F_n \right) \left( \tmop{Id} +
      \sum_{n \geq 1} G_n \right) & \\
      & = &\di \tmop{Id} + \sum_{n \geq 1} \sum_{k = 0}^n F_k G_{n - k} & ( F_0
      = G_0 = \tmop{Id} )
    \end{array}
  \end{equation}
  But for $l \geq 1$, $G_l ( x ) = \psi_l x^{l + 1}$ and then, for $k \geq 1$,
  \begin{equation}
    \begin{array}{ccc}
      F_k G_l .x & = & \di \sum_{\tmmathbf{n}= ( n_1, \ldots, n_s ) \in
      \mathcal{N}_k} \frac{1}{l ( \tmmathbf{n} ) !} \varphi_{\tmmathbf{n}}
      x^{\| \tmmathbf{n} \| + l ( \tmmathbf{n} )} \partial_x^{l ( \tmmathbf{n}
      )} ( \psi_l x^{l + 1} )\\
      & = &\di  \sum_{\tmmathbf{n}= ( n_1, \ldots, n_s ) \in \mathcal{N}_k}
      \psi_l \frac{1}{l ( \tmmathbf{n} ) !} \varphi_{\tmmathbf{n}} \frac{( l +
      1 ) !}{( l + 1 - l ( \tmmathbf{n} ) ) !} x^{\| \tmmathbf{n} \| + l +
      1}\\
      & = &\di  \left( \sum_{\tmmathbf{n}= ( n_1, \ldots, n_s ) \in
      \mathcal{N}_k} C^{l ( \tmmathbf{n} )}_{l + 1} \varphi_{\tmmathbf{n}}
      \psi_l \right) x^{k + l + 1}
    \end{array}
  \end{equation}
  and then, for $n \geq 1$,
  \begin{equation}
    \eta_n = \varphi_n + \psi_n + \sum_{k = 1}^{n - 1} \sum_{\tmmathbf{n}= (
    n_1, \ldots, n_s ) \in \mathcal{N}_k} C^{l ( \tmmathbf{n} )}_{l + 1}
    \varphi_{\tmmathbf{n}} \psi_{n - k}
  \end{equation}
  If now $\tilde{\varphi} = \varphi^{- 1}$ and $\tilde{F} = \tau^{- 1} (
  \tilde{\varphi} )$, then, as $\tilde{F} F = \tmop{Id}$ we get
  \begin{equation}
    \tilde{F} = \tmop{Id} + \sum_{s \geq 1} ( - 1 )^s F_{n_1} \ldots F_{n_s} =
    \tmop{Id} + \sum_{\tmmathbf{n} \in \mathcal{N}} ( - 1 )^{l ( \tmmathbf{n}
    )} F_{\tmmathbf{n}}
  \end{equation}
  but for $\tmmathbf{n}= ( n_1, \ldots, n_s ) \in \mathcal{N}_{}$,
  \begin{equation}
    F_{\tmmathbf{n}} ( x ) = \sum_{\tmscript{\begin{array}{c}
      \tmmathbf{m}^i \in \mathcal{N}_{n_i}\\
      1\leq i\leq s\\
      \end{array}}} B_1 ( \tmmathbf{m}^1, \ldots,^{} \tmmathbf{m}^s )
    \varphi_{\tmmathbf{m}^1} \ldots \varphi_{\tmmathbf{m}^s} x^{\|
    \tmmathbf{n} \| + 1}
  \end{equation}
  Now
  \begin{equation}
    \tilde{\varphi}_n = \sum_{\text{$\tmmathbf{n}= ( n_1, \ldots, n_s ) \in
    \mathcal{N}_n$}} ( - 1 )^s \sum_{\tmscript{\begin{array}{c}
      \tmmathbf{m}^i \in \mathcal{N}_{n_i}\\
     1\leq i\leq s\\
      \end{array}}} B_1 ( \tmmathbf{m}^1, \ldots,^{} \tmmathbf{m}^s )
    \varphi_{\tmmathbf{m}^1} \ldots \varphi_{\tmmathbf{m}^s} \label{e38}
  \end{equation}
  and this gives the attempted result.

\end{proof}

Using the same ideas, we will finally prove theorems \ref{th2} and \ref{th3}

\subsection{Proof of Theorems \ref{th2} and \ref{th3}}

As before, let $\varphi ( x ) = x + \di \sum_{n \geq 1} \varphi_n x^{n + 1}$ and
$\psi ( x ) = x + \di \sum_{n \geq 1} \psi_n x^{n + 1}$ two elements of $G_2$ and
$\eta = \mu ( \varphi, \psi ) = \psi \circ \varphi$ with
\begin{equation}
  \eta ( x ) = x + \sum_{n \geq 1} \eta_n x^{n + 1}
\end{equation}
If
\begin{equation}
  \begin{array}{cccccc}
    f ( x ) & = & \log ( \varphi' ( x ) ) & = & \di \sum_{n \geq 1}
    \frac{f_n}{n!} 
    x^n & \hspace{1em} ( f_n = \gamma_n ( \varphi ) )\\
    g ( x ) & = & \log ( \psi' ( x ) ) & = & \di \sum_{n \geq 1} \frac{g_n}{n!}
    x^n & \hspace{1em} ( g_n = \gamma_n ( \psi ) )\\
    h ( x ) & = & \log ( \eta' ( x ) ) & = & \di \sum_{n \geq 1} \frac{h_n}{n!}
    x^n & \hspace{1em} ( h_n = \gamma_n ( \eta ) )
  \end{array}
\end{equation}
then
\begin{equation}
  \begin{array}{rcl}
    h ( x ) & = & \log ( ( \psi \circ \varphi )' ( x ) )\\
    & = & \log ( \varphi' ( x ) . \psi' ( \varphi ( x ) )\\
    & = & \log ( \varphi' ( x ) ) + ( \log \psi' ) \circ \varphi ( x )\\
    & = & f ( x ) + F ( g ) ( x )
  \end{array}
\end{equation}
where $F$ is the substitution automorphism associated to $\varphi$. We remind
that $F = \tmop{Id} + \sum_{n \geq 1} F_n$. Because of equation \ref{e19},
\begin{equation}
  \begin{array}{ccc}
    F_n & = & \di \sum_{\tmmathbf{n}= ( n_1, \ldots, n_s ) \in \mathcal{N}_n}
    \frac{1}{l ( \tmmathbf{n} ) !} \varphi_{\tmmathbf{n}} x^{\| \tmmathbf{n}
    \| + l ( \tmmathbf{n} )} \partial_x^{l ( \tmmathbf{n} )}\\
    & = & \di \sum_{\tmmathbf{n}= ( n_1, \ldots, n_s ) \in \mathcal{N}_n}
    \frac{1}{l ( \tmmathbf{n} ) !} \sum_{\tmscript{\begin{array}{c}
      \tmmathbf{m}^i \in \mathcal{N}_{n_i}\\
      1\leq i\leq s\\
      \end{array}}} \frac{A ( \tmmathbf{m}^1, \ldots, \tmmathbf{m}^s )
    f_{\tmmathbf{m}^1} \ldots f_{\tmmathbf{m}^s}}{ \tmmathbf{m}^1 ! \ldots
    \tmmathbf{m}^s !} x^{\| \tmmathbf{n} \| + l ( \tmmathbf{n} )}
    \partial_x^{l ( \tmmathbf{n} )}\\
    & = & \di \sum_{\tmmathbf{n} \in \mathcal{N}_n}
    \frac{f_{\tmmathbf{n}}}{\tmmathbf{n} !} \sum_{\tmmathbf{m}^1 \ldots
    \tmmathbf{m}^s = \tmmathbf{n}} A ( \tmmathbf{m}^1, \ldots, \tmmathbf{m}^s
    ) \frac{1}{s!} x^{n + s} \partial_x^s
  \end{array} \label{e42}
\end{equation}
But for $k \geq 1$,
\begin{equation}
  \begin{array}{ccc}
    F_n \left( \frac{g_k}{k!} x^k \right) & = & 
  \end{array} \sum_{\tmmathbf{n} \in \mathcal{N}_n} \frac{f_{\tmmathbf{n}}
  g_k}{\tmmathbf{n} !k!} \sum_{\tmmathbf{m}^1 \ldots \tmmathbf{m}^s =
  \tmmathbf{n}} A ( \tmmathbf{m}^1, \ldots, \tmmathbf{m}^s ) C_k^s x^{n + k}
\end{equation}
and we obtain immediately the formula for the coproduct.

Let now $\tilde{\varphi} = \varphi^{- 1}$ and
\begin{equation}
  \tilde{f} ( x ) = 
    \log ( \tilde{\varphi}' ( x ) )  = \sum_{n \geq 1}
    \frac{\tilde{f}_n}{n!} x^n  \quad ( \tilde{f}_n = \gamma_n (
    \tilde{\varphi} ) )
\end{equation}
Since $\tilde{\varphi} \circ \varphi ( x ) = x$,
\begin{equation}
  0 = \log ( ( \tilde{\varphi} \circ \varphi )' ( x ) ) = f ( x ) + F.
  \tilde{f} ( x )
\end{equation}
thus
\begin{equation}
  \tilde{f} ( x ) = - \tilde{F} .f ( x ) = - f ( x ) - \sum_{\tmmathbf{n} \in
  \mathcal{N}} ( - 1 )^{l ( \tmmathbf{n} )} F_{\tmmathbf{n}} ( f ) ( x )
\end{equation}
But, once again,
\begin{equation}
  \begin{array}{lll}
    f_{\tmmathbf{n}, k} ( x ) & = & \di \sum_{\tmmathbf{n} \in \mathcal{N}_n} ( -
    1 )^{l ( \tmmathbf{n} )} F_{\tmmathbf{n}} ( \frac{f_k}{k!} x^k )\\
    & = & \di \sum_{\tmmathbf{n} \in \mathcal{N}_n} ( - 1 )^{l ( \tmmathbf{n} )}
    \frac{f_k}{k!} \sum_{\tmscript{\begin{array}{c}
      \tmmathbf{m}^i \in \mathcal{N}_{n_i}\\
      1\leq i\leq s\\
      \end{array}}} B_k ( \tmmathbf{m}^1, \ldots,^{} \tmmathbf{m}^s )
    \varphi_{\tmmathbf{m}^1} \ldots \varphi_{\tmmathbf{m}^s} x^{\|
    \tmmathbf{n} \| + k}\\
    & = & \di \sum_{\tmmathbf{n} \in \mathcal{N}_n} \sum_{\tmmathbf{m}^1 \ldots
    \tmmathbf{m}^s = \tmmathbf{n}} ( - 1 )^s B_k ( \tmmathbf{m}^1, \ldots,^{}
    \tmmathbf{m}^s ) \varphi_{\tmmathbf{n}} \frac{f_k}{k!} x^{\| \tmmathbf{n}
    \| + k}\\
    & = & - \di \sum_{\tmmathbf{n} \in \mathcal{N}_n} U_k ( \tmmathbf{n} )
    \varphi_{\tmmathbf{n}} \frac{f_k}{k!} x^{\| \tmmathbf{n} \| + k}
  \end{array}
\end{equation}
Now, replacing $\varphi_{\tmmathbf{n}}$ as in equation \ref{e42},
\begin{equation}
  \begin{array}{lll}
    f_{\tmmathbf{n}, k} ( x ) & = & \di \sum_{\tmmathbf{n} \in \mathcal{N}_n} ( -
    1 )^{l ( \tmmathbf{n} )} F_{\tmmathbf{n}} ( \frac{f_k}{k!} x^k )\\
    & = & - \di \sum_{\tmmathbf{n} \in \mathcal{N}_n} U_k ( \tmmathbf{n} )
    \varphi_{\tmmathbf{n}} \frac{f_k}{k!} x^{\| \tmmathbf{n} \| + k}\\
    & = & - \di \sum_{\tmmathbf{n} \in \mathcal{N}_n} \frac{f_{\tmmathbf{n}}
    f_k}{\tmmathbf{n} !k!} \sum_{\tmmathbf{m}^1 \ldots \tmmathbf{m}^s =
    \tmmathbf{n}} A ( \tmmathbf{m}^1, \ldots, \tmmathbf{m}^s ) U_k ( \|
    \tmmathbf{m}^1 \|, \ldots, || \tmmathbf{m}^s \| ) x^{n + k}
  \end{array}
\end{equation}
Now, for $l \geq 1$,
$$
     \tilde{f}_l  =  - f_l + \sum_{n = 1}^{l - 1} \sum_{\tmmathbf{n} \in
     \mathcal{N}_n} \frac{ l!f_{\tmmathbf{n}} f_{l - n}}{\tmmathbf{n} ! ( l - n )
     !} \sum_{\tmmathbf{m}^1 \ldots \tmmathbf{m}^s = \tmmathbf{n}} A (
     \tmmathbf{m}^1, \ldots, \tmmathbf{m}^s ) U_{l - n} ( \| \tmmathbf{m}^1
     \|, \ldots, || \tmmathbf{m}^s \| )
$$
and this gives immediately the attempted formula.

This ends the proofs for our formulas but, as we said before, the first proofs
were derived from mould calculus and we will give the main ideas in the next
sections.

\section{Mould calculus and the shuffle Hopf algebra $\tmop{sh}
(\mathbbm{N}^{\ast} )$.\label{moulds}}

\subsection{An example of mould calculus}

\subsubsection{Formal Conjugacy of equations}

Mould calculus, as defined by J. Ecalle (see
{\cite{et2}},{\cite{et3}},{\cite{esn}}), appears in the study of formal or
analytic conjugacy of differential equations, vector fields, diffeomorphisms.
In order to introduce it, we give here a very simple but useful example.

Let $u \in G_2$ and the associated equation
\[ ( E_u ) \hspace{2em} \partial_t x = u ( x ) = x + \sum_{n \geq 1} u_n x^{n
   + 1} \]
For $u$ and $v$ in $G_2$ the equations $( E_u )$ and $( E_v )$ are formally
conjugated if there exists an element $\varphi$ of $G_2$ such that, if $x$ is
a solution of $( E_u )$ then $y = \varphi ( x )$ is a solution of $( E_v )$.
This defines an equivalence relation on the set of such equations and one can
easily check that there is only one class : For any equation $( E_u )$, there
exist a unique $\varphi_{}$ of $G_2$ such that, if $x$ is a solution of $( E_u
)$ then $y = \varphi_{} ( x )$ is a solution of
\[ ( E_0 ) \hspace{2em} \partial_t y = y \]
The equation for $\varphi_{}$ reads
\begin{equation}
  u ( x ) \varphi'_{} ( x ) = \varphi_{} ( x )
\end{equation}
and, if
\begin{equation}
  \varphi ( x ) = x + \sum_{n \geq 1} \varphi_n x^{n + 1}
\end{equation}
then
\begin{equation}
  \begin{array}{ccc}
    u_1 + 2 \varphi_1 & = & \varphi_1\\
    u_2 + 2 \varphi_1 u_1 + 3 \varphi_2 & = & \varphi_2\\
    & \vdots & \\
    u_n + \sum_{k = 1}^{n - 1} ( k + 1 ) u_{n - k} \varphi_k + ( n + 1 )
    \varphi_n & = & \varphi_n
  \end{array}
\end{equation}
Recursively, one can determine the values $a_n ( \varphi ) = \varphi_n$ and
thus the diffeomorphism $\varphi$. This does not give a direct formula for the
coefficients of $\varphi$. Among other properties that may be useful for more
sophisticated equations, we will see that the mould calculus will give
explicit formulas.

Mould calculus, for this example, is based on two remarks which are detailed
in the next two sections.

\subsubsection{Diffeomorphisms an substitution automorphisms}

As we have seen in section \ref{aut}, to any diffeomorphism $\varphi \in G_2$
one can associate a substitution automorphism $F \in \tilde{G}_2$
\begin{equation}
  F = \tmop{Id} + \sum_{n \geq 1} F_n
\end{equation}
Moreover, the action of such an operator on a product of formal power series
induces a coproduct
\begin{equation}
  \Delta F_{} = F_{} \otimes F_{} \hspace{1em} ( F_{} ( f g ) = ( F_{} f ) (
  F_{} g ) )
\end{equation}
which also reads
\begin{equation}
  \forall n \geq 1, \hspace{1em} \Delta F_n = F_n \otimes \tmop{Id} + \sum_{k
  = 1}^{n - 1} F_k \otimes F_{n - k} + \tmop{Id} \otimes F_n
\end{equation}

\subsubsection{Symmetral moulds and shuffle Hopf algebra}

Now, for $u \in G_2$, the equation $( E_u )$ reads
\begin{equation}
  \partial_t x = \left( \mathbbm{B}_0 + \sum_{n \geq 1} u_n \mathbbm{B}_n
  \right) .x =\mathbbm{B}.x \hspace{1em} \text{with} \hspace{1em}
  \mathbbm{B}_n = x^{n + 1} \partial_x
\end{equation}
Instead of computing the conjugating map $\varphi$ we could look for its
associated substitution automorphism $F_{}$ in the following shape :
\begin{equation}
  F = \tmop{Id} + \sum_{s \geq 1} \sum_{n_1 \geq 1, \ldots, n_s \geq 1}
  M^{n_1, \ldots, n_s} \mathbbm{B}_{n_1} \ldots \mathbbm{B}_{n_s}
\end{equation}
As we will see later, in order to get a substitution automorphism, is is
sufficient to impose that for any sequences $\tmmathbf{k}= ( k_1, \ldots, k_s
)$ and $\tmmathbf{l}= ( l_1, \ldots, l_t )$,
\begin{equation}
  M^{\tmmathbf{k}} M^{\tmmathbf{l}} = \sum_{\tmmathbf{m}}
  \tmop{sh}^{\tmmathbf{k},\tmmathbf{l} }_{\tmmathbf{m}} M^{\tmmathbf{m}}
\end{equation}
where $\tmop{sh}^{\tmmathbf{k},\tmmathbf{l} }_{\tmmathbf{m}}$ is the
number of shuffling of the sequences 
$\tmmathbf{k},\tmmathbf{l}$ that gives the sequence $\tmmathbf{m}$. The set of
such coefficients is called a symmetral mould. Moreover the conjugacy equation
reads
\begin{equation}
  \mathbbm{B}F_{} .x = F_{} \mathbbm{B}_0 .x
\end{equation}
Now we can solve the equation $\mathbbm{B}F_{} = F_{} \mathbbm{B}_0$
by noticing
that, for $( n_1, \ldots, n_s ) \in (\mathbbm{N}^{\ast} )^s$,
\begin{equation}
  \left[ \mathbbm{B}_0,\mathbbm{B}_{n_1} \ldots \mathbbm{B}_{n_s} \right] = (
  n_1 + \ldots + n_s )\mathbbm{B}_{n_1} \ldots \mathbbm{B}_{n_s}
\end{equation}
and using this commutation relations, one can check that for $s = 1$ and a
sequence $( n_1$) we get
\begin{equation}
  u_{n_1} + n_1 M^{n_1} = 0
\end{equation}
and for $s \geq 2$ and a sequence $( n_1, \ldots, n_s ) \in
(\mathbbm{N}^{\ast} )^s$,
\begin{equation}
  u_{n_1} M^{n_2, \ldots, n_s} + ( n_1 + \ldots + n_s ) M^{n_1, \ldots, n_s} =
  0
\end{equation}
This defines a symmetral mould, for $s \geq 1$ and $( n_1, \ldots, n_s ) \in
(\mathbbm{N}^{\ast} )^s$,
\begin{equation}
  M^{n_1, \ldots, n_s} = \frac{( - 1 )^s u_{n_1} \ldots u_{n_s}}{( n_1 +
  \ldots + n_s ) ( n_2 + \ldots + n_s ) \ldots ( n_{s - 1} + n_s ) n_s}
\end{equation}
thus we get explicit formulas for $F$ and $\varphi ( x ) = F_{} .x$ : For $n
\geq 1$,
\begin{equation}
  \varphi_n x^{n + 1} = \sum_{s = 1}^n \sum_{\tmscript{\begin{array}{c}
    n_1 + \ldots + n_s = n\\
    n_i \geq 1
  \end{array}}} M^{n_1, \ldots, n_s} \mathbbm{B}_{n_1} \ldots
  \mathbbm{B}_{n_s} .x
\end{equation}
and
\begin{equation}
  \varphi_n = \sum_{s = 1}^n \sum_{\tmscript{\begin{array}{c}
    n_1 + \ldots + n_s = n\\
    n_i \geq 1
  \end{array}}} ( n_s + 1 ) ( n_{s - 1} + n_s + 1 ) \ldots ( n_2 + \ldots +
  n_s + 1 ) M^{n_1, \ldots, n_s}
\end{equation}

We just gave the outlines of the method here. The important idea is that we
only used the commutation of $\mathbbm{B}_0$ with the over derivations
$\mathbbm{B}_n$ ($n \geq 1$), which means that we worked as these derivations
were free of other relations. This can be interpreted in the following
algebraic way.

\subsection{The free group and its Hopf algebra of coordinates}

\subsubsection{Lie algebra and substitution automorphisms}

Let $\mathcal{A}^1$ the Lie algebra of formal vector fields generated by the
derivations

\begin{equation}
  \forall n \geq 1, \quad \mathbbm{B}_n = x^{n + 1} \partial_x
\end{equation}
Its enveloping algebra $\mathcal{U}(\mathcal{A}^1 )$ is a graded Hopf algebra
and, see {\cite{cm}}, the Hopf algebra $\mathcal{H}^1_{\tmop{CM}}$ is the dual
of $\mathcal{U}(\mathcal{A}^1 )$. Note that this dual is well-defined as the
graded components of $\mathcal{U}(\mathcal{A}^1 )$ are vector spaces
of finite dimension. If $G (\mathcal{A}^1 ) \subset
\mathcal{U}(\mathcal{A}^1 )$  is
the group of the group-like elements of $\mathcal{U}(\mathcal{A}^1 )$, this is
exactly the group of substitution automorphism describe above and it is
isomorphic to the group $G_2$
\begin{equation}
  \forall F \in G (\mathcal{A}^1 ), \forall f \in \mathbbm{R}[ [ x ] ] \quad
  F.f = f \circ \varphi^{}, \quad \varphi \in G_2
\end{equation}
In other terms, $G (\mathcal{A}^1 ) = \widetilde{G_{}}_2$.

\subsubsection{The free group and its Hopf algebra of coordinates}

Our previous mould calculus suggests to introduce, by analogy with
$\mathcal{A}^1$, the graded free Lie algebra $A^1$ generated by a set of
primitive elements $X_n$, $n \geq 1$,
\begin{equation}
  \Delta ( X_n ) = X_n \otimes 1 + 1 \otimes X_n
\end{equation}
The enveloping algebra $\mathcal{U}( A^1 )$ is a Hopf algebra which is also
called the concatenation Hopf algebra in combinatorics (see {\cite{reut}}). If
the unity is $X_{\emptyset} = 1$ ($\emptyset$ is the empty sequence), then an
element $\tmmathbf{U}$ of $\mathcal{U}( A^1 )$ can be written
\begin{equation}
  \begin{array}{lll}
    \tmmathbf{U} & = & \di U^{\emptyset} X_{\emptyset} + \sum_{s \geq 1}
    \sum_{n_1, \ldots, n_s \geq 1} U^{n_1, \ldots, n_s} \mathcal{} X_{n_1}
    \ldots X_{n_s}\\
    & = &\di  U^{\emptyset} X_{\emptyset} + \sum_{s \geq 1} \sum_{n_1, \ldots,
    n_s \geq 1} U^{n_1, \ldots, n_s} \mathcal{} X_{n_1, \ldots, n_s}\\
    & = &\di  \sum U^{\bullet} X_{\bullet}
  \end{array}
\end{equation}
where the collection of coefficients $U^{\bullet}$ is called a
\textit{mould.} The structure of the enveloping algebra $\mathcal{U}( A^1 )$
can be described as follows : the product is given by
\begin{equation}
  \forall \tmmathbf{m},\tmmathbf{n} \in \mathcal{N}, \quad X_{\tmmathbf{m}}
  X_{\tmmathbf{n}} = X_{\tmmathbf{m} \tmmathbf{n}} \quad
  \text{(concatenation),}
\end{equation}
the coproduct is
\begin{equation}
  \Delta ( X_{\tmmathbf{n}} ) = \sum_{\tmmathbf{n}^1,\tmmathbf{n}^2} \tmop{sh}
  \left( \begin{array}{c}
    \tmmathbf{n}^1,\tmmathbf{n}^2\\
    \tmmathbf{n}
  \end{array} \right) X_{\tmmathbf{n}^1} \otimes X_{\tmmathbf{n}^2}
\end{equation}
where $\tmop{sh}^{\tmmathbf{n}^1,\tmmathbf{n}^2}_{\tmmathbf{n}}$ is
the number of shuffling of the sequences 
$\tmmathbf{n}^1,\tmmathbf{n}^2$ that gives $\tmmathbf{n}$. Finally, the
antipode $S$ is defined by
\begin{equation}
  S ( X_{n_1, \ldots, n_s} ) = ( - 1 )^s X_{n_s, \ldots, n_1}
\end{equation}

Once again one can define the group $G ( A^1 )$ and if $\tmmathbf{F} \in G (
A^1 )$ then
\begin{equation}
  \tmmathbf{F}= \sum_{\tmmathbf{n} \in \mathcal{N} \cup \{ \emptyset \}}
  F^{\tmmathbf{n}} X_{\tmmathbf{n}}
\end{equation}
where the mould $F^{\bullet}$ is \textit{symmetral} : $F^{\emptyset} = 1$ and
\begin{equation}
  \forall \tmmathbf{n}^1,\tmmathbf{n}^2, \quad F^{\tmmathbf{n}^1}
  F^{\tmmathbf{n}^2} = \sum_{\tmmathbf{n}}
  \tmop{sh}^{\tmmathbf{n}^1,\tmmathbf{n}^2}_{\tmmathbf{n}} F^{\tmmathbf{n}}
\end{equation}
Moreover, if $\tmmathbf{G}$ is the group inverse of $\tmmathbf{F}$, then its
associated mould is given by the formulas
\[ G^{n_1, \ldots, n_s} = ( - 1 )^s F^{n_s, \ldots, n_1} \]
Thanks to the graduation on $\mathcal{U}( A^1 )$, its dual $H^1$ is a Hopf
algebra, the Hopf algebra of coordinates on $G ( A^1 )$ and, if the dual basis
of $\{ X_{\tmmathbf{n}}, \hspace{1em} \tmmathbf{n} \in \mathcal{N} \}$ is $\{
Z^{\tmmathbf{n}}, \hspace{1em} \tmmathbf{n} \in \mathcal{N} \}$ then the
product in $H^1$ is defined by :
\begin{equation}
  \forall \tmmathbf{n}^1,\tmmathbf{n}^2, \quad Z^{\tmmathbf{n}^1}
  Z^{\tmmathbf{n}^2} = \sum_{\tmmathbf{n}}
  \tmop{sh}^{\tmmathbf{n}^1,\tmmathbf{n}^2}_{\tmmathbf{n} }
   Z^{\tmmathbf{n}}
\end{equation}
The coproduct is :
\begin{equation}
  \Delta ( Z^{\tmmathbf{n}} ) = Z^{\tmmathbf{n}} \otimes 1 + 1 \otimes
  Z^{\tmmathbf{n}} + \sum_{\tmmathbf{n}^1 \tmmathbf{n}^2 =\tmmathbf{n}}
  Z^{\tmmathbf{n}^1} \otimes Z^{\tmmathbf{n}^2}
\end{equation}
where $\tmmathbf{n}^1 \tmmathbf{n}^2$ is the concatenation of the two nonempty
sequences $\tmmathbf{n}^1$ and $\tmmathbf{n}^2$ and $Z^{\emptyset} = 1$ is the
unity. Finally, the antipode is given by
\begin{equation}
  S ( Z^{n_1, \ldots, n_s} ) = ( - 1 )^s Z^{n_s, \ldots, n_1}
\end{equation}
The structure of $H^1$ (coproduct, antipode, ...) is fully explicit. This will
be of great use since our previous mould calculus suggests that there exists a
surjective morphism from $A^1$ on $\mathcal{A}^1$ that induces an injective
morphism from $\mathcal{H}^1_{\tmop{CM}}$ into $H^1$. In other words,
$\mathcal{H}^1_{\tmop{CM}}$ can be identified to a sub-Hopf algebra of $H^1$
and, as everything is explicit in $H^1$, one can derive formulas for the
coproduct and the antipode in $\mathcal{H}_{\tmop{CM}}^1$.

\section{Morphisms.\label{mor}}

The application defined by $\rho ( X_n ) =\mathbbm{B}_n = x^{n + 1}
\partial_x$ obviously determines a morphism from $A^1$ (resp. $\mathcal{U}(
A^1 )$, resp. $G ( A^1 )$) on $\mathcal{A}^1$ (resp.
$\mathcal{U}(\mathcal{A}^1 )$, resp. $G (\mathcal{A}^1 ) \simeq G_2$) and it is
surjective : If $\varphi \in G_2$ and $F = \tau^{- 1} ( \varphi ) \in G
(\mathcal{A}^1 ) = \tilde{G}_2$, then, if
\begin{equation}
  b ( x ) = x + \sum_{n \geq 1} b_n x^{n + 1} = \frac{\varphi ( x )}{\varphi'
  ( x )}
\end{equation}
then $\varphi$ is the unique diffeomorphism of $G_2$ that conjugates $( E_b )$
to $( E_0 )$ thus
\begin{equation}
  F = \tmop{Id} + \sum_{\tmmathbf{n} \in \mathcal{N}} M^{\tmmathbf{n}}
  \mathbbm{B}_{\tmmathbf{n}} = \rho \left( X_{\emptyset} + \sum_{\tmmathbf{n}
  \in \mathcal{N}} M^{\tmmathbf{n}} X_{\tmmathbf{n}} \right)
\end{equation}
By duality, it induces a morphism $\rho^{\ast}$ from  $\mathcal{H}^1$ to $H^1$
by
\begin{equation}
  \forall \gamma \in \mathcal{H}^1, \quad \rho^{\ast} ( \gamma ) = \gamma
  \circ \rho
\end{equation}
and, since $\rho$ is surjective, $\rho^{\ast}$ is injective :
$\mathcal{H}_{\tmop{CM}}^1$ is isomorphic to the sub-Hopf algebra $\rho^{\ast}
(\mathcal{H}^1_{\tmop{CM}} ) \subset H^1$. Using this injective morphism, we
define
\begin{equation}
  \forall n \geq 1, \quad \Gamma_n = \rho^{\ast} ( \gamma_n )
\end{equation}
and $\rho^{\ast} (\mathcal{H}_{\tmop{CM}}^1 )$ is then the Hopf algebra
generated by the $\Gamma_n$. In order to get formulas in
$\mathcal{H}_{\tmop{CM}}^1$, we will use the algebra $\rho^{\ast}
(\mathcal{H}^1_{\tmop{CM}} )$ and express the $\Gamma_n$ in terms of the
$Z^{\tmmathbf{n}}$ :

\begin{theorem}
  \label{th8}For $n \geq 1$,
  \begin{equation}
    \begin{array}{ccc}
      \Gamma_n & = & \di n! \sum_{\tmscript{\begin{array}{c}
        \tmmathbf{n}= ( n_1, \ldots, n_s ) \in \mathcal{N}_n
      \end{array}}} \sum_{t = 1}^s \frac{( - 1 )^{t - 1}}{t}
      \sum_{\tmscript{\begin{array}{c}
        \tmmathbf{n}^1 \ldots \tmmathbf{n}^t =\tmmathbf{n}
      \end{array}}} Z^{\tmmathbf{n}^1} \ldots Z^{\tmmathbf{n}^t}
      S^{\tmmathbf{n}^1} \ldots S^{\tmmathbf{n}^t}\\
      & = &\di  n! \sum_{\tmscript{\begin{array}{c}
        \tmmathbf{n}= ( n_1, \ldots, n_s ) \in \mathcal{N}_n
      \end{array}}} Q^{\tmmathbf{n}} Z^{\tmmathbf{n}}
    \end{array}
  \end{equation}
  where $S^{n_1, \ldots, n_s} = \prod_{i = 1}^s ( n_i + n_{i + 1} + \ldots +
  n_s + 1 ) = \prod_{i = 1}^s ( \hat{n}_i + 1 )$ and $Q^{n_1, \ldots, n_s} = (
  n_s + 1 ) \prod_{i = 2}^s \hat{n}_i$ with $Q^{n_1} = ( n_1 + 1 )$.
\end{theorem}

Let $\tmmathbf{F}= X_{\emptyset} + \sum_{\tmmathbf{n} \in \mathcal{N}}
F^{\tmmathbf{n}} X_{\tmmathbf{n}} \in G ( A^1 )$. If $F = \rho (\tmmathbf{F})
\in G (\mathcal{A}^1 )$, then
\begin{equation}
  \Gamma_n (\tmmathbf{F}) = \gamma_n ( F ) = \gamma_n ( \varphi ) = (
  \partial^n_x \log ( \varphi' ) ( x ) )_{x = 0}
\end{equation}
where $\varphi \in G_2$ is defined by :
\begin{equation}
  \varphi^{} ( x ) = \rho (\tmmathbf{F}) .x = F.x = x +
  \sum_{\tmscript{\begin{array}{c}
    ( n_1, \ldots, n_s ) \in \mathcal{N}
  \end{array}}} F^{n_1, \ldots, n_s} \mathbbm{B}_{n_1} \ldots
  \mathbbm{B}_{n_s} .x
\end{equation}
Then
\begin{equation}
  \varphi' ( x ) = 1 + \sum_{\tmscript{\begin{array}{c}
    ( n_1, \ldots, n_s ) \in \mathcal{N}
  \end{array}}} F^{n_1, \ldots, n_s} S^{n_1, \ldots, n_s} x^{n_1 + \ldots +
  n_s}
\end{equation}
Using the logarithm and derivation, one easily gets the formula
\begin{equation}
  \Gamma_n (\tmmathbf{F}) = n! \sum_{\tmscript{\begin{array}{c}
    \tmmathbf{n}= ( n_1, \ldots, n_s ) \in \mathcal{N}_n
  \end{array}}} \sum_{t = 1}^s \frac{( - 1 )^{t - 1}}{t}
  \sum_{\tmscript{\begin{array}{c}
    \tmmathbf{n}^1 \ldots \tmmathbf{n}^t =\tmmathbf{n}
  \end{array}}} F^{\tmmathbf{n}^1} \ldots F^{\tmmathbf{n}^t}
  S^{\tmmathbf{n}^1} \ldots S^{\tmmathbf{n}^t}
\end{equation}

We prove the second part of the formula in section \ref{ini}, using the fact
that $F^{\bullet}$ is symmetral. As
\begin{equation}
  \Gamma_n (\tmmathbf{F}) = n! \sum_{\tmscript{\begin{array}{c}
    \tmmathbf{n}= ( n_1, \ldots, n_s ) \in \mathcal{N}_n
  \end{array}}} F^{\tmmathbf{n}} Q^{\tmmathbf{n}}
\end{equation}
and $Z^{\tmmathbf{n}} .\tmmathbf{F}= F^{\tmmathbf{n}}$, theorem \ref{th8} will
be proved.

As $\rho^{\ast} ( \delta_n ) = \Gamma_n \in H^1$, and, since the coproduct and
the antipode are explicit in $H^1$, we can once again obtain the formulas
given in theorems \ref{th2} and \ref{th3}.

\section{Initial Proofs.} \label{ini}

\subsection{Proof of theorem \ref{th8}}

We already proved that, for $n \geq 1$,
\begin{equation}
  \Gamma_n = n! \sum_{\tmscript{\begin{array}{c}
    \tmmathbf{n}= ( n_1, \ldots, n_s ) \in \mathcal{N}_n
  \end{array}}} \sum_{t = 1}^s \frac{( - 1 )^{t - 1}}{t}
  \sum_{\tmscript{\begin{array}{c}
    \tmmathbf{n}^1 \ldots \tmmathbf{n}^t =\tmmathbf{n}
  \end{array}}} Z^{\tmmathbf{n}^1} \ldots Z^{\tmmathbf{n}^t}
  S^{\tmmathbf{n}^1} \ldots S^{\tmmathbf{n}^t}
\end{equation}
Extending the notion of shuffling, for $t \geq 1$, if $\tmmathbf{m}^1,
\ldots,\tmmathbf{m}^t,\tmmathbf{m}$ are $t + 1$ sequences, then $\tmop{sh}^{
  \tmmathbf{m}^1, \ldots,\tmmathbf{m}^t}_{\tmmathbf{m}}$ is the number
of ways to obtain the sequence 
$\tmmathbf{m}$ by shuffling the sequences $\tmmathbf{m}^1,
\ldots,\tmmathbf{m}^t$. Then,
\begin{equation}
  \begin{array}{ccc}
   \di  \frac{1}{n!} \Gamma_n & = &\di  \sum_{\tmmathbf{n} \in
     \mathcal{N}_n} \sum_{t 
    = 1}^{l (\tmmathbf{n})} \frac{( - 1 )^{t - 1}}{t}
    \sum_{
      \tmmathbf{n}^1 \ldots \tmmathbf{n}^t =\tmmathbf{n}}
    Z^{\tmmathbf{n}^1} \ldots Z^{\tmmathbf{n}^t} 
    S^{\tmmathbf{n}^1} \ldots S^{\tmmathbf{n}^t}\\
    & = &\di  \sum_{\tmmathbf{n} \in \mathcal{N}_n} \sum_{t = 1}^{l
    (\tmmathbf{n})} \frac{( - 1 )^{t - 1}}{t} \sum_{
      \tmmathbf{n}^1 \ldots \tmmathbf{n}^t =\tmmathbf{n}} \left(
      \sum_{\tmmathbf{m}}  \tmop{sh}^{\tmmathbf{n}^1,
        \ldots,\tmmathbf{n}^t}_{\tmmathbf{m} }
   Z^{\tmmathbf{m}} \right) S^{\tmmathbf{n}^1} \ldots
    S^{\tmmathbf{n}^t}\\
    & = &\di  \sum_{\tmmathbf{m} \in \mathcal{N}_n} \left( Z^{\tmmathbf{m}}
    \sum_{t = 1}^{l (\tmmathbf{m})} \frac{( - 1 )^{t - 1}}{t}
    \sum_{\tmmathbf{n}^1, \ldots,\tmmathbf{n}^t \in \mathcal{N}}
    \tmop{sh}^{\tmmathbf{n}^1, \ldots,\tmmathbf{n}^t }_{
      \tmmathbf{m}}
    S^{\tmmathbf{n}^1} \ldots S^{\tmmathbf{n}^t} \right.
  \end{array}
\end{equation}
Note that in these equations, we had $\| \tmmathbf{m} \| = \| \tmmathbf{n}
\|$ and $l (\tmmathbf{m}) = l (\tmmathbf{n})$. For a given sequence
$\tmmathbf{m} \in \mathcal{N}$, let
\begin{equation}
  Q^{\tmmathbf{m}} = \sum_{t = 1}^{l (\tmmathbf{m})} \frac{( - 1 )^{t - 1}}{t}
  \sum_{\tmmathbf{n}^1, \ldots,\tmmathbf{n}^t \in \mathcal{N}}
  \tmop{sh}^{\tmmathbf{n}^1, \ldots,\tmmathbf{n}^t }_{
    \tmmathbf{m}}
   S^{\tmmathbf{n}^1} \ldots S^{\tmmathbf{n}^t}
\end{equation}
it remains to prove that, if $\tmmathbf{m}= ( m_1, \ldots, m_s )$ then
$Q^{m_1, \ldots, m_s} = ( m_s + 1 ) \prod_{i = 2}^s \hat{m}_i$ with $Q^{m_1} =
( m_1 + 1 )$. We prove this formula by induction on $l (\tmmathbf{m})$.

If $l (\tmmathbf{m}) = 1$, then $\tmmathbf{m}= ( m_1 )$ and
\begin{equation}
  Q^{m_1} = \frac{( - 1 )^0}{1} \sum_{\tmmathbf{n}^1 \in \mathcal{N}}
  \tmop{sh}^{ \tmmathbf{n}^1}_{
    \tmmathbf{m}}
   S^{\tmmathbf{n}^1} = S^{m_1} = m_1 + 1
\end{equation}
If $l (\tmmathbf{m}) = s \geq 2$, then let $\tmmathbf{m}= ( m_1, \ldots, m_s
)$ and $\tmmathbf{p}= ( m_2, \ldots, m_s )$. For any sequence $\tmmathbf{n}= (
n_1, \ldots, n_k )$, we note $m_1 \tmmathbf{n}= ( m_1, n_1, \ldots, n_k )$. If
a shuffling of $t \geq 1$ sequences $\tmmathbf{n}^1, \ldots,\tmmathbf{n}^t$
gives $\tmmathbf{m}$ then
\begin{itemize}
  \item Either there exists $1 \leq i \leq t$ such that $\tmmathbf{n}^i = (
  m_1 )$ (but then $t \geq 2$), and, omitting $\tmmathbf{n}^i= (
  m_1 )$, the
  corresponding shuffling of the $t - 1$ remaining sequences gives
  $\tmmathbf{p}$.
  
  \item Either there exists $1 \leq i \leq t$ such that $\tmmathbf{n}^i = m_1
  \tilde{\tmmathbf{n}}^i$ ($\tilde{\tmmathbf{n}}^i \not= \emptyset$)
  (necessarily, $t < l (\tmmathbf{m})$) and, replacing $\tmmathbf{n}^i$ by
  $\tilde{\tmmathbf{n}}^i$, the corresponding shuffling of the $t$ sequences
  gives $\tmmathbf{p}$.
\end{itemize}
This means that :
\begin{equation}
  \begin{array}{rcl}
    Q^{\tmmathbf{m}} & = & \di \sum_{t = 1}^{l (\tmmathbf{m})} \frac{( - 1 )^{t -
    1}}{t} \sum_{\tmmathbf{n}^1, \ldots,\tmmathbf{n}^t}
\tmop{sh}^{\tmmathbf{n}^1, \ldots,\tmmathbf{n}^t }_{
      \tmmathbf{m}}
     S^{\tmmathbf{n}^1} \ldots S^{\tmmathbf{n}^t}\\
    & = & \di \sum_{t = 1}^{l (\tmmathbf{m}) - 1} \frac{( - 1 )^{t - 1}}{t}
    \sum_{\tmmathbf{n}^1, \ldots,\tmmathbf{n}^t}
    \tmop{sh}^{\tmmathbf{n}^1, \ldots,\tmmathbf{n}^t }_{
      \tmmathbf{p}}
     \sum_{i = 1}^t S^{\tmmathbf{n}^1} \ldots S^{m_1
    \tmmathbf{n}^i} \ldots S^{\tmmathbf{n}^t}\\
    &  & +\di  \sum_{t = 2}^{l (\tmmathbf{m})} \frac{( - 1 )^{t - 1}}{t}
    \sum_{\tmmathbf{n}^1, \ldots,\tmmathbf{n}^{t - 1}}
    \tmop{sh}^{\tmmathbf{n}^1, \ldots,\tmmathbf{n}^{t - 1} }_{\tmmathbf{p}}
     \sum_{i = 0}^{t - 1} S^{\tmmathbf{n}^1} \ldots
    S^{\tmmathbf{n}^i} S^{m_1} S^{\tmmathbf{n}^{i + 1}} \ldots
    S^{\tmmathbf{n}^{t - 1}}
  \end{array}
\end{equation}
but as $S^{m_1} = m_1 + 1$ and $S^{m_1 \tmmathbf{n}^i} = ( m_1 + \|
\tmmathbf{n}^i \| + 1 ) S^{\tmmathbf{n}^i}$,
\begin{equation}
  \begin{array}{rcl}
    Q^{\tmmathbf{m}} & = & \di \sum_{t = 1}^{l (\tmmathbf{m}) - 1} \frac{( - 1
    )^{t - 1}}{t} \sum_{\tmmathbf{n}^1, \ldots,\tmmathbf{n}^t}
  \tmop{sh}^{\tmmathbf{n}^1, \ldots,\tmmathbf{n}^t }_{
      \tmmathbf{p}}
     \sum_{i = 1}^t ( m_1 + \| \tmmathbf{n}^i \| + 1 )
    S^{\tmmathbf{n}^1} \ldots S^{\tmmathbf{n}^i} \ldots S^{\tmmathbf{n}^t}\\
    &  & + \di \sum_{t = 1}^{l (\tmmathbf{m}) - 1} \frac{( - 1 )^t}{t + 1}
    \sum_{\tmmathbf{n}^1, \ldots,\tmmathbf{n}^t}
    \tmop{sh}^{\tmmathbf{n}^1, \ldots,\tmmathbf{n}^t }_{
      \tmmathbf{p}}
     \sum_{i = 0}^t ( m_1 + 1 ) S^{\tmmathbf{n}^1} \ldots
    S^{\tmmathbf{n}^i} S^{\tmmathbf{n}^{i + 1}} \ldots S^{\tmmathbf{n}^t}\\
    & = & \di \sum_{t = 1}^{l (\tmmathbf{m}) - 1} \frac{( - 1 )^{t - 1}}{t} ( t (
    m_1 + 1 ) + \| \tmmathbf{p} \| ) \sum_{\tmmathbf{n}^1,
    \ldots,\tmmathbf{n}^t} \tmop{sh}^{\tmmathbf{n}^1,
    \ldots,\tmmathbf{n}^t }_{
      \tmmathbf{p}}
    S^{\tmmathbf{n}^1} \ldots S^{\tmmathbf{n}^t}\\
    &  & + \di \sum_{t = 1}^{l (\tmmathbf{m}) - 1} \frac{( - 1 )^t}{t + 1} ( t +
    1 ) ( m_1 + 1 ) \sum_{\tmmathbf{n}^1, \ldots,\tmmathbf{n}^t}
    \tmop{sh}^{\tmmathbf{n}^1, \ldots,\tmmathbf{n}^t }_{
      \tmmathbf{p}}
     S^{\tmmathbf{n}^1} \ldots S^{\tmmathbf{n}^t}\\
    & = & \di \| \tmmathbf{p} \| \sum_{t = 1}^{l (\tmmathbf{p})} \frac{( - 1 )^{t
    - 1}}{t} \sum_{\tmmathbf{n}^1, \ldots,\tmmathbf{n}^t}
\tmop{sh}^{\tmmathbf{n}^1, \ldots,\tmmathbf{n}^t }_{
      \tmmathbf{p}}
    S^{\tmmathbf{n}^1} \ldots S^{\tmmathbf{n}^t}\\
    & = & \| \tmmathbf{p} \| Q^{\tmmathbf{p}}
  \end{array}
\end{equation}
And it obviously gives the right formula for $Q^{\tmmathbf{m}}$.

\subsection{Proof of theorem \ref{th2}}

Using the above formula we have
\begin{equation}
  \begin{array}{rcl}
    \Delta \Gamma_n & = & \di n! \sum_{\tmmathbf{m} \in \mathcal{N}_n}
    Q^{\tmmathbf{m}}_{} ( \Delta Z^{\tmmathbf{m}} )\\
    & = & \di n! \sum_{\tmmathbf{m} \in \mathcal{N}_n} Q^{\tmmathbf{m}}_{} \left(
    Z^{\tmmathbf{m}} \otimes 1 + 1 \otimes Z^{\tmmathbf{m}} +
    \sum_{\tmscript{\begin{array}{c}
      \tmmathbf{p} \tmmathbf{q}=\tmmathbf{m}
    \end{array}}} Z^{\tmmathbf{p}} \otimes Z^{\tmmathbf{q}} \right)\\
    & = & \di \left( n! \sum_{\tmmathbf{m} \in \mathcal{N}_n} Q^{\tmmathbf{m}}_{}
    Z^{\tmmathbf{m}} \right) \otimes 1 + 1 \otimes \left( n!
    \sum_{\tmmathbf{m} \in \mathcal{N}_n} Q^{\tmmathbf{m}}_{} Z^{\tmmathbf{m}}
    \right)\\
    &  & \di + n! \sum_{\tmmathbf{m} \in \mathcal{N}_n}
    \sum_{\tmscript{\begin{array}{c}
      \tmmathbf{p} \tmmathbf{q}=\tmmathbf{m}
    \end{array}}} Q^{\tmmathbf{m}}_{} Z^{\tmmathbf{p}} \otimes
    Z^{\tmmathbf{q}}\\
    & = & \di \Gamma_n \otimes 1 + 1 \otimes \Gamma_n + n! \sum_{\tmmathbf{m} \in
    \mathcal{N}_n} \sum_{\tmscript{\begin{array}{c}
      \tmmathbf{p} \tmmathbf{q}=\tmmathbf{m}
    \end{array}}} Q^{\tmmathbf{m}}_{} Z^{\tmmathbf{p}} \otimes
    Z^{\tmmathbf{q}}\\
    & = & \di \Gamma_n \otimes 1 + 1 \otimes \Gamma_n + \tilde{\Delta} \Gamma_n
  \end{array}
\end{equation}
Now if $\tmmathbf{p} \tmmathbf{q}=\tmmathbf{m}= ( m_1, \ldots, m_s )$ with
$\tmmathbf{p},\tmmathbf{q} \in \mathcal{N}$ ($s \geq 2$), then
\begin{equation}
  \frac{Q^{\tmmathbf{m}}}{Q^{\tmmathbf{q}}} = \frac{( m_s + 1 ) \prod_{i =
  2}^s \hat{m}_i}{( m_s + 1 ) \prod^s_{i = l (\tmmathbf{p}) + 2} \hat{m}_i} =
  \prod_{i = 2}^{l (\tmmathbf{p}) + 1} \hat{m}_i = \prod_{i = 2}^{l
  (\tmmathbf{p}) + 1} ( \hat{p}_i + \| \tmmathbf{q} \| ) =
  R^{\tmmathbf{p}}_{\| \tmmathbf{q} \|}
\end{equation}
with the convention that if $i = l (\tmmathbf{p}) + 1$, then $\hat{p}_i = 0$.
As this coefficient only depends on $\tmmathbf{p}$ and $\| \tmmathbf{q} \|$,
\begin{equation}
  \begin{array}{rcl}
    \tilde{\Delta} \Gamma_n & = &\di  n! \sum_{\tmmathbf{m} \in \mathcal{N}_n}
    \sum_{\tmscript{\begin{array}{c}
      \tmmathbf{p} \tmmathbf{q}=\tmmathbf{m}
    \end{array}}} Q^{\tmmathbf{m}}_{} Z^{\tmmathbf{p}} \otimes
    Z^{\tmmathbf{q}}\\
    & = &\di  n! \sum_{\tmmathbf{m} \in \mathcal{N}_n}
    \sum_{\tmscript{\begin{array}{c}
      \tmmathbf{p} \tmmathbf{q}=\tmmathbf{m}
    \end{array}}} ( R^{\tmmathbf{p}}_{\| \tmmathbf{q} \|} Z^{\tmmathbf{p}} )
    \otimes Q^{\tmmathbf{q}} Z^{\tmmathbf{q}}\\
    & = &\di  n! \sum_{k = 1}^{n - 1} \left( \sum_{\tmmathbf{p} \in
    \mathcal{N}_{n - k}} R^{\tmmathbf{p}}_{\| \tmmathbf{q} \|}
    Z^{\tmmathbf{p}} \right) \otimes \left( \sum_{\tmmathbf{q} \in
    \mathcal{N}_k} Q^{\tmmathbf{q}} Z^{\tmmathbf{q}} \right)\\
    & = &\di  \sum_{k = 1}^{n - 1} \left( \frac{n!}{k!} \sum_{\tmmathbf{p} \in
    \mathcal{N}_k} R^{\tmmathbf{p}}_{\| \tmmathbf{q} \|} Z^{\tmmathbf{p}}
    \right) \otimes \Gamma_k\\
    & = &\di  \sum_{k = 1}^{n - 1} P^n_k \otimes \Gamma_k
  \end{array}
\end{equation}
and it remains to prove that, for $n \geq 1$ and $1 \leq k \leq n - 1$,
\begin{equation}
  P^n_k = \sum_{\tmscript{\begin{array}{c}
    ( n_1, \ldots, n_s ) \in \mathcal{N}\\
    n_1 + \ldots + n_s = n - k, s \geq 1
  \end{array}}} \frac{n!}{n_1 ! \ldots n_s !k!} \alpha^{n_1, \ldots, n_s}_k
  \Gamma_{n_1} \ldots \Gamma_{n_s}
\end{equation}
with
\begin{equation}
  \alpha^{\tmmathbf{n}}_k = \sum_{t = 1}^{l (\tmmathbf{n})} C^{t_{}}_k
  \sum_{\tmscript{\begin{array}{c}
    \tmmathbf{n}^1 \ldots \tmmathbf{n}^t =\tmmathbf{n}\\
    \tmmathbf{n}^i \not= \emptyset
  \end{array}}} \frac{1}{l (\tmmathbf{n}^1 ) ! \ldots l (\tmmathbf{n}^t ) !}
  \prod_{i = 1}^t \frac{1}{\| \tmmathbf{n}^i \| + 1}
\end{equation}
This formula was first conjectured on the first values of $n$. Now let
\begin{equation}
  \begin{array}{rcl}
    \tilde{P}^n_k & = & \di \sum_{\tmscript{\begin{array}{c}
      ( n_1, \ldots, n_s ) \in \mathcal{N}_{n - k}
    \end{array}}} \frac{n!}{n_1 ! \ldots n_s !k!} \alpha^{n_1, \ldots, n_s}_k
    \Gamma_{n_1} \ldots \Gamma_{n_s}\\
    & = & \di \frac{n!}{k!} \sum_{\tmscript{\begin{array}{c}
      ( n_1, \ldots, n_s ) \in \mathcal{N}_{n - k}
    \end{array}}} \alpha^{n_1, \ldots, n_s}_k \sum_{\tmmathbf{m}^i \in
    \mathcal{N}_{n_i}} Q^{\tmmathbf{m}^1} \ldots Q^{\tmmathbf{m}^s}
    Z^{\tmmathbf{m}^1} \ldots Z^{\tmmathbf{m}^s}\\
    & = & \di \frac{n!}{k!} \sum_{\tmscript{\begin{array}{c}
      ( n_1, \ldots, n_s ) \in \mathcal{N}_{n - k}
    \end{array}}} \alpha^{n_1, \ldots, n_s}_k \sum_{\tmmathbf{p}}
    \sum_{\tmmathbf{m}^i \in \mathcal{N}_{n_i}}
    \tmop{sh}^{\tmmathbf{m}^1, \ldots,\tmmathbf{m}^s }_{
      \tmmathbf{p}}
     Q^{\tmmathbf{m}^1} \ldots Q^{\tmmathbf{m}^s}
    Z^{\tmmathbf{p}}\\
    & = & \di \frac{n!}{k!} \sum_{\tmmathbf{p} \in \mathcal{N}_{n - k}}
    Z^{\tmmathbf{p}} \sum_{s \geq 1} \sum_{\tmmathbf{m}^1,
    \ldots,\tmmathbf{m}^s} \alpha^{\| \tmmathbf{m}^1 \|, \ldots, \|
    \tmmathbf{m}^s \|}_k \tmop{sh}^{\tmmathbf{m}^1,
    \ldots,\tmmathbf{m}^s }_{
      \tmmathbf{p}}
     Q^{\tmmathbf{m}^1} \ldots Q^{\tmmathbf{m}^s}
  \end{array}
\end{equation}
It remains to prove that for a given $\tmmathbf{p} \in \mathcal{N}_{n - k}$,
we have
\[ \tilde{R}^{\tmmathbf{p}}_k = \sum^{l (\tmmathbf{p})}_{s = 1}
   \sum_{\tmmathbf{m}^1, \ldots,\tmmathbf{m}^s} \alpha^{\| \tmmathbf{m}^1 \|,
   \ldots, \| \tmmathbf{m}^s \|}_k \tmop{sh}^{\tmmathbf{m}^1,
   \ldots,\tmmathbf{m}^s }_{
     \tmmathbf{p}}
    Q^{\tmmathbf{m}^1} \ldots Q^{\tmmathbf{m}^s} =
   R^{\tmmathbf{p}}_k = \prod_{i = 2}^{l (\tmmathbf{p}) + 1} ( \hat{p}_i + k )
\]
As in the previous proof, if $l (\tmmathbf{p}) = 1$ then $R^{p_1}_k = k$ and
\begin{equation}
  \tilde{R}^{p_1}_k = \alpha^{p_1}_k Q^{p_1} = C^1_k \frac{1}{l (\tmmathbf{p})
  !} \frac{1}{p_1 + 1} ( p_1 + 1 ) = k
\end{equation}
and if $l (\tmmathbf{p}) \geq 2$, as $\tmmathbf{p}= p_1 \tmmathbf{q}$,
\begin{equation}
  \begin{array}{rcl}
    \tilde{R}^{\tmmathbf{p}}_k & = & \tilde{R}^{p_1 \tmmathbf{q}}_k\\
    & = & \di \sum^{l (\tmmathbf{p}) - 1}_{s = 1} \sum_{\tmmathbf{m}^1,
    \ldots,\tmmathbf{m}^s \atop 1\leq i\leq s}  \alpha^{\| \tmmathbf{m}^1 \|,
    \ldots, \| \tmmathbf{m}^i \| + p_1, \ldots, \| \tmmathbf{m}^s \|}_k
    \tmop{sh}^{\tmmathbf{m}^1, \ldots,\tmmathbf{m}^s }_{
      \tmmathbf{q}}
    Q^{\tmmathbf{m}^1} \ldots Q^{p_1 \tmmathbf{m}^i} \ldots
    Q^{\tmmathbf{m}^s}\\
    &  & + \di \sum^{l (\tmmathbf{p}) - 1}_{s = 1}
    \hspace{-2mm}\sum_{\tmmathbf{m}^1, 
     \ldots,\tmmathbf{m}^s\atop 0\leq i\leq s} \hspace{-2mm}
   \alpha^{\| \tmmathbf{m}^1 \|, 
    \ldots, \| \tmmathbf{m}^i \|, p_1, \| \tmmathbf{m}^{i + 1} \|, \ldots, \|
    \tmmathbf{m}^s \|}_k \tmop{sh}^{\tmmathbf{m}^1, \ldots,\tmmathbf{m}^s}_{
      \tmmathbf{q}}
    Q^{\tmmathbf{m}^1} \ldots Q^{\tmmathbf{m}^s} Q^{p_1}
  \end{array}
\end{equation}
Since $Q^{p_1} = ( p_1 + 1 )$ and $Q^{p_1 \tmmathbf{m}^i} = \| \tmmathbf{m}^i
\| Q^{\tmmathbf{m}^i}$, we get
\begin{equation}
  \tilde{R}^{\tmmathbf{p}}_k = \sum^{l (\tmmathbf{p}) - 1}_{s = 1}
  \sum_{\tmmathbf{m}^1, \ldots,\tmmathbf{m}^s}
  \tmop{sh}^{\tmmathbf{m}^1, \ldots,\tmmathbf{m}^s }_{
    \tmmathbf{q}}
   Q^{\tmmathbf{m}^1} \ldots Q^{\tmmathbf{m}^s} V^{\|
  \tmmathbf{m}^1 \|, \ldots, \| \tmmathbf{m}^s ||}_{k, p_1}
\end{equation}
where
\begin{equation}
  V^{n_1, \ldots, n_s}_{k, p_1} = \sum_{i = 1}^s n_i\alpha^{n_1, \ldots, n_i +
  p_1 \ldots, n_s}_k  + ( p_1 + 1 ) \sum_{i = 0}^s \alpha^{n_1, \ldots,
  n_i, p_1, n_{i + 1} \ldots, n_s}_k
\end{equation}
but
\begin{equation}
  \begin{array}{rcl}
    V^{n_1, \ldots, n_s}_{k, p_1} & = & \di \sum_{i = 1}^s n_i \sum_{t = 1}^s
    C^{t_{}}_k \sum_{\tmscript{\begin{array}{c}
      \tmmathbf{n}^1 \ldots \tmmathbf{n}^t = ( n_1, \ldots, n_i + p_1 \ldots,
      n_s )
    \end{array}}} A (\tmmathbf{n}^1, \ldots,\tmmathbf{n}^t )\\
    &  & \di+ ( p_1 + 1 ) \sum_{i = 0}^s \sum_{t = 1}^{s + 1} C^{t_{}}_k
    \sum_{\tmscript{\begin{array}{c}
      \tmmathbf{n}^1 \ldots \tmmathbf{n}^t = ( n_1, \ldots, n_i, p_1, \ldots,
      n_s )
    \end{array}}} A (\tmmathbf{n}^1, \ldots,\tmmathbf{n}^t )
  \end{array}
\end{equation}
In the first term, we get a sequence $\tmmathbf{n}^1 \ldots \tmmathbf{n}^t = (
n_1, \ldots, n_i + p_1 \ldots, n_s )$ starting with a decomposition
$\tmmathbf{m}^1 \ldots \tmmathbf{m}^t = ( n_1, \ldots, n_s )$ and adding $p_1$
to one element of one of the sequences $\tmmathbf{m}^i$. In the second term,
$\tmmathbf{n}^1 \ldots \tmmathbf{n}^t = ( n_1, \ldots, n_i, p_1, \ldots,
n_s )$, then either $p_1$ is one of the sequences $\tmmathbf{n}^1
\ldots \tmmathbf{n}^t$, and, once it is omitted, we
get a decomposition $\tmmathbf{m}^1 \ldots \tmmathbf{m}^{t - 1} = ( n_1,
\ldots, n_s )$, either we start with a decomposition $\tmmathbf{m}^1 \ldots
\tmmathbf{m}^t = ( n_1, \ldots, n_s )$ and $p_1$ is inserted in one of the
sequences $\tmmathbf{m}^i$ : If $\tmmathbf{n}=( n_1, \ldots, n_s )$, then,
\begin{equation}
  \begin{array}{rcl}
    V^{n_1, \ldots, n_s}_{k, p_1} & = &\di \sum_{t \geq 1}^{} C^{t_{}}_k
    \sum_{\tmscript{\begin{array}{c}
      \tmmathbf{n}^1 \ldots \tmmathbf{n}^t = \tmmathbf{n}
    \end{array}}} \sum_{i = 1}^t \frac{\| \tmmathbf{n}^i \| ( \|
    \tmmathbf{n}^i \| + 1 )}{\| \tmmathbf{n}^i \| + p_1 + 1} A
    (\tmmathbf{n}^1, \ldots,\tmmathbf{n}^t )\\
    &  &\di  + ( p_1 + 1 ) \sum_{t \geq 1}^{} C^{t + 1_{}}_k
    \sum_{\tmscript{\begin{array}{c}
      \tmmathbf{n}^1 \ldots \tmmathbf{n}^t = \tmmathbf{n}
    \end{array}}} \frac{t + 1}{p_1 + 1} A (\tmmathbf{n}^1,
    \ldots,\tmmathbf{n}^t )\\
    &  &\di  + ( p_1 + 1 ) \sum_{t \geq 1}^{} C^{t_{}}_k
    \sum_{\tmscript{\begin{array}{c}
      \tmmathbf{n}^1 \ldots \tmmathbf{n}^t = \tmmathbf{n}
    \end{array}}} \sum_{i = 1}^t \frac{\| \tmmathbf{n}^i \| + 1}{\|
    \tmmathbf{n}^i \| + p_1 + 1} A (\tmmathbf{n}^1, \ldots,\tmmathbf{n}^t )
  \end{array}
\end{equation}
But
\begin{equation}
  ( p_1 + 1 ) C^{t + 1_{}}_k \frac{t + 1}{p_1 + 1} = C^t_k ( k - t )
\end{equation}
and
\begin{equation}
  \frac{\| \tmmathbf{n}^i \| ( \| \tmmathbf{n}^i \| + 1 )}{\| \tmmathbf{n}^i
  \| + p_1 + 1} + ( p_1 + 1 ) \frac{\| \tmmathbf{n}^i \| + 1}{\|
  \tmmathbf{n}^i \| + p_1 + 1} = \| \tmmathbf{n}^i \| + 1
\end{equation}
thus
\begin{equation}
  \begin{array}{rcl}
    V^{n_1, \ldots, n_s}_{k, p_1} & = & \di \sum_{t \geq 1}^{} C^{t_{}}_k
    \sum_{
      \tmmathbf{n}^1 \ldots \tmmathbf{n}^t = \tmmathbf{n}
    } A (\tmmathbf{n}^1, \ldots,\tmmathbf{n}^t ) \left( ( k - t )
    + \sum_{i = 1}^t ( \| \tmmathbf{n}^i \| + 1 ) \right)  \\
    & = & \di ( n_1 + \ldots + n_s + k ) \alpha^{n_1, \ldots, n_s}_k  
  \end{array}
\end{equation}
Now by induction we get, if $\tmmathbf{p}= p_1 \tmmathbf{q}$,
\begin{equation}
  \begin{array}{rcl}
    \tilde{R}^{\tmmathbf{p}}_k & = & \di \sum^{l (\tmmathbf{p}) - 1}_{s = 1}
    \sum_{\tmmathbf{m}^1, \ldots,\tmmathbf{m}^s}
    \tmop{sh}^{\tmmathbf{m}^1, \ldots,\tmmathbf{m}^s }_{
      \tmmathbf{q}}
     Q^{\tmmathbf{m}^1} \ldots Q^{\tmmathbf{m}^s} ( \|
    \tmmathbf{q} \| + k ) \alpha^{\| \tmmathbf{m}^1 \|, \ldots, \|
    \tmmathbf{m}^s \|}_k\\
    & = &\di  ( \| \tmmathbf{q} \| + k ) \tilde{R}^{\tmmathbf{q}}_k\\
    & = &\di  ( \| \tmmathbf{q} \| + k ) R^{\tmmathbf{q}}_k\\
    & = &\di  ( \| \tmmathbf{q} \| + k ) \prod_{i = 2}^{l (\tmmathbf{q}) + 1} (
    \hat{q}_i + k )\\
    & = &\di  \prod_{i = 1}^{l (\tmmathbf{q}) + 1} ( \hat{q}_i + k )\\
    & = &\di  \prod_{i = 2}^{l (\tmmathbf{p}) + 1} ( \hat{q}_i + k )\\
    & = &\di  R^{\tmmathbf{p}}_k
  \end{array}
\end{equation}

We live the second proof of theorem \ref{th3} to the reader : the ideas are
the same, noticing that
\[ \begin{array}{rcl}
     S ( \Gamma_n ) & = & \di \sum_{( n_1, \ldots, n_s ) \in \mathcal{N}_n}
     Q^{n_1, \ldots, n_s} S ( Z^{n_1, \ldots, n_s} )\\
     & = & \di \sum_{( n_1, \ldots, n_s ) \in \mathcal{N}_n} ( - 1 )^s Q^{n_1,
     \ldots, n_s} Z^{n_s, \ldots, n_1}
   \end{array} \]

\section{Tables and conclusion.}

Some computations give the following tables.

\subsection{The coproduct}

The table gives the value of $\frac{n!}{n_1 ! \ldots n_{s + 1} !} \alpha^{n_1,
\ldots, n_s}_{n_{s + 1}}$ for a given sequence $( n_1, \ldots, n_{s + 1} )$ :
\begin{table}[!h]
\begin{tabular}{|c|c|c|c|}
  \hline
  $( 1, 1 ) = 1$ &  &  &  \\
  \hline
  $( 1, 2 ) = 3$ & $( 2, 1 ) = 1$ & $( 1, 1, 1 ) = 1$ &   \\
  \hline
  $( 1, 3 ) = 6$ & $( 2, 2 ) = 4$ & $( 3, 1 ) = 1$ & $( 1, 1, 2 ) = 7$
  \\
\hline
 $( 1,2, 1 ) = 3 / 2$&$( 2, 1, 1 ) = 3 / 2$ & $( 1, 1, 1, 1 ) = 1$ &\\
  \hline
  $( 1, 4 ) = 10$ & $( 2, 3 ) = 10$ & $( 3, 2 ) = 5$ & $( 4, 1 ) = 1$
  \\
\hline
 $( 1,1, 3 ) = 25$ & $( 1, 3, 1 ) = 2$ & $( 3, 1, 1 ) = 2$ & $( 1, 2,
 2 ) = 25 / 2$ \\
\hline 
 $( 2, 1, 2 ) = 25 / 2$ & $( 2, 2, 1 ) = 3$& $( 1, 1, 1, 2 ) = 15$ &
 $( 1, 1, 2, 1 ) = 2$ \\
\hline
 $( 1, 2, 1, 1 ) = 2$ & $( 2,
  1, 1, 1 ) = 2$ & $( 1, 1, 1, 1, 1 ) = 1$ &\\
  \hline
\end{tabular}
\end{table}

This gives
\[ \begin{array}{ccc}
     \tilde{\Delta} \Gamma_1 & = & 0\\
     \tilde{\Delta} \Gamma_2 & = &  \Gamma_1
     \otimes \Gamma_1 \\
     \tilde{\Delta} \Gamma_3 & = &  ( \Gamma_2 + \Gamma_1^2 )
     \otimes \Gamma_1 + 3 \Gamma_1 \otimes \Gamma_2 \\
     \tilde{\Delta} \Gamma_4 & = &   ( \Gamma_3 + 3 \Gamma_1^{}
     \Gamma_2 + \Gamma_1^3 ) \otimes \Gamma_1 + ( 4 \Gamma_2 + 7 \Gamma_1^2 )
     \otimes \Gamma_2 + 6 \Gamma_1 \otimes \Gamma_3 \\
     \tilde{\Delta} \Gamma_5 & = & ( \Gamma_4 + 4
     \Gamma_1
     \Gamma_3 + 3 \Gamma_2^2 + 6 \Gamma^2_1 \Gamma_2 + \Gamma_1^4 ) \otimes
     \Gamma_1\\
     &  & + ( 5 \Gamma_3 + 25 \Gamma_1^{} \Gamma_2 + 15 \Gamma_1^3 ) \otimes
     \Gamma_2 + ( 10 \Gamma_2 + 25 \Gamma_1^2 ) \otimes \Gamma_3 + 10 \Gamma_1
     \otimes \Gamma_4 
   \end{array} \]

\subsection{The antipode}

The table gives the value of $\frac{( n_1 + \ldots + n_s ) !}{n_1 ! \ldots n_s
!} \beta^{n_1, \ldots, n_s}_{}$ for a given sequence $( n_1, \ldots, n_s )$ :
\begin{table}[!h]
\begin{tabular}{|c|c|c|c|}
  \hline
  $( 1 ) = - 1$ &  &  &   \\
  \hline
  $( 2 ) = - 1$ & $( 1, 1 ) = 1$ &  &  \\
  \hline
  $( 3 ) = - 1$ & $( 1, 2 ) = 3$ & $( 2, 1 ) = 1$ & $( 1, 1, 1 ) = - 2$  \\
  \hline
  $( 4 ) = - 1$ & $( 1, 3 ) = 6$ & $( 2, 2 ) = 4$ & $( 3, 1 ) = 1$  \\
  \hline
  $( 1, 1, 2 ) = - 11$ & $( 1, 2, 1 ) = - 9 / 2$ & $( 2, 1, 1 ) = - 5 / 2$ &
  $( 1, 1, 1, 1 ) = 6$  \\
  \hline
  $( 5 ) = - 1$ & $( 1, 4 ) = 10$ & $( 2, 3 ) = 10$ & $( 3, 2 ) = 5$  \\
  \hline
  $( 4, 1) = 1$& $( 1, 1, 3 ) = - 35$ & $( 1, 3, 1 ) = - 8$ & $( 3, 1,
  1 ) = - 3$ \\ 
  \hline
  $( 1, 2, 2 ) = - 55 / 2$ & $( 2, 1, 2 ) = - 35 / 2$ & $( 2, 2, 1 ) = - 7$ & 
   \\
  \hline
  $( 1, 1, 1, 2 ) = 50$ & $( 1, 1, 2, 1 ) = 22$ & $( 1, 2, 1, 1 ) = 29 / 2$ &
  $( 2, 1, 1, 1 ) = 19 / 2$  \\
  \hline
  $( 1, 1, 1, 1, 1 ) = - 24$ &  &  &   \\
  \hline
\end{tabular}
\end{table}

This gives :
\[ \begin{array}{lll}
     S ( \Gamma_1 ) & = & - \Gamma_1\\
     S ( \Gamma_2 ) & = & - \Gamma_2 + \Gamma_1^2\\
     S ( \Gamma_3 ) & = & - \Gamma_3 + 4 \Gamma_1 \Gamma_2 - 2 \Gamma_1^3\\
     S ( \Gamma_4 ) & = & - \Gamma_4 + 7 \Gamma_1 \Gamma_3 + 4 \Gamma_2^2 - 18
     \Gamma_1^2 \Gamma_2 + 6 \Gamma_1^4\\
     S ( \Gamma_5 ) & = & - \Gamma_5 + 11 \Gamma_1 \Gamma_4 + 15 \Gamma_2
     \Gamma_3 - 46 \Gamma_1^2 \Gamma_3 - 52 \Gamma_1 \Gamma_2^2 + 96
     \Gamma_1^3 \Gamma_2 - 24 \Gamma_1^5
   \end{array} \]

This is the attempted result but the formulas in proposition \ref{prop7},
theorem \ref{th2} and \ref{th3} are not unique because
$\mathcal{H}^1_{\tmop{CM}}$ is commutative and, in the computations, it is
much more ''simple'' to consider that the algebra generated by the $\delta_n$
is somehow noncommutative. This situation calls for furthers investigations,
since the coefficients appearing in proposition \ref{prop7} for the
Faà di Bruno coordinates seem to arise in the study of a
noncommutative version of diffeomorphisms (see {\cite{fra}}).

\bibliographystyle{plain}
\bibliography{main}

\end{document}